\documentclass[a4paper,english,oneside]{article}
\usepackage[T1]{fontenc}
\usepackage[utf8]{inputenc} 
\usepackage{authblk}

\usepackage{geometry}
\usepackage{amsmath}
\usepackage{amsfonts}
\usepackage{amssymb}
\usepackage{amsthm}
\usepackage[english]{babel}
\usepackage{color}
\usepackage{lmodern}
\usepackage{pstricks-add}
\usepackage{hyperref}
\hypersetup{hidelinks}

\usepackage{mathscinet}
\usepackage{enumitem}
\usepackage[cal=boondoxo,bb=ams]{mathalfa}
\usepackage{microtype}
\usepackage{mathtools}
\usepackage{esint}

\newtheorem{theorem}{Theorem}[section]

\newtheorem{definition}{Definition}[section]
\newtheorem{lemma}{Lemma}[section]
\newtheorem{proposition}{Proposition}[section]
\newtheorem{remark}{Remark}[section]

\numberwithin{equation}{section}

\begin{document}
\title{On the blow-up of solutions to a Nakao-type problem with a time-dependent damping term}

\author[a,b]{Yuequn Li}
\author[b]{Alessandro Palmieri\footnote{
\href{mailto:alessandro.palmieri@uniba.it}{alessandro.palmieri@uniba.it}}}

\affil[a]{School of Mathematical Sciences, Nanjing Normal University Nanjing 210023, China}
\affil[b]{Department of Mathematics, University of Bari, Via E. Orabona 4, 70125 Bari, Italy}

\date{}

\maketitle

\begin{abstract}
In this paper, we study a semilinear weakly coupled system of wave equations with power nonlinearities. More precisely, we couple (through the nonlinear terms) a wave equation and a damped wave equation with a time-dependent coefficient for the damping term. For the coefficient of the damping term we consider two cases:  the scale-invariant case and the scattering producing case. By applying an iteration argument, we get a blow-up result and  upper bound estimates for the lifespan of the solutions. In the scale-invariant case, we obtain a shift of the space dimension in the blow-up region for the same weakly coupled system with a classical damping (i.e. with a constant coefficient), while for the scattering producing case we find the same blow-up region as for the classical Nakao problem.
\end{abstract}


\begin{flushleft}
\textbf{Keywords} weakly coupled system, blow-up, lifespan estimates, iteration argument.
\end{flushleft}

\begin{flushleft}
\textbf{AMS Classification (2020)} 35B44, 35G55, 35L52, 35L71.
\end{flushleft}

\section{Introduction}\label{s1}

In this paper, we consider the blow-up of local solutions to the following weakly coupled system of semilinear wave equations
\begin{equation}\label{eqs}
\begin{cases}
\partial_t^2u-\Delta u+\partial_tu=\vert v\vert^p, &  t>0, \  x\in \mathbb{R}^n,\\
\partial_t^2v-\Delta v+b(t)\partial_tv=\vert u\vert^q, &  t>0, \  x\in \mathbb{R}^n,\\
(u,\partial_tu,v,\partial_tv)(0,x)=\varepsilon(u_0,u_1,v_0,v_1)(x),  &x\in\mathbb{R}^n,
\end{cases}
\end{equation}
where $p,q>1$ and $\varepsilon>0$ describes the size of the initial data. 

Our analysis will be carried out separately for the following two cases:
\begin{enumerate}
\item a scale-invariant coefficient for the damping term $b(t)=\frac{\mu}{1+t}$, with $\mu>0$;
\item a scattering producing coefficient for the damping term $b \in L^1([0,\infty))$, $b\geq 0$.
\end{enumerate}
Let us give a short overview on these two categories of cofficients for the damping term in the $v$-equation.

\noindent (1) The influence of a scale-invariant damping term on the properties of the solutions to (linear and semilinear) wave equations have been extensively studied over the last two decades. For the Cauchy problem associated to the wave equation $\partial_t^2 u-\Delta u+\frac{\mu}{1+t}\partial_t u=0$, known in the literature also as Euler-Poisson-Darboux equation (up to a translation in the space variable), Wirth \cite{Wi04} established an explicit representation formula by using the theory of special functions and the phase space analysis and obtained energy and $L^p-L^{p'}$ estimates. Furthermore,  in \cite{Pa2021MMAS,YG2021JDE} integral representation formulas are established for the solution of the inhomogeneous linear wave equation with a scale-invariant damping term by using Yagdjian's integral transform approach. For semilinear wave equations with a scale-invariant damping term and power nonlinearity, several blow-up and global existence results have been obtained, see \cite{Ab2015,Da2021,AbLuRe2015,IkSo2018,LaiPaTa2025,LaiTaWa2017,Pa2024,Wa2014} and references therein. Roughly speaking, the coefficient $\frac{\mu}{1+t}$ has a bordeline nature among the time-dependent coefficients for the damping term of a wave equation (cf. \cite{Wi04,Wi06,Wi07}). Thus, the constant $\mu$ plays a crucial role in the decay rates for the $L^p-L^{p'}$ estimates of the linear model and, consequently, in the value of the critical exponent for the power nonlinearity $|u|^p$. By critical exponent, since we are dealing with nonnegative power nonlinearities, we mean the threshold value that separates the blow-up range from the small data solutions' global existence range.

\noindent (2) Let us recall briefly the reason why in the case $b\in L^1([0,\infty))$ we say that we are considering a scattering producing damping term $b(t)\partial_t v$. In \cite[Theorem 3.1]{Wiphdt} it is proved the following result: \emph{if $b=b(t)$ is a nonnegative and summable function, then there exists an isomorphism $\mathcal{W}:\mathcal{E}\to \mathcal{E}$, where $\mathcal{E}\doteq \dot{H}^1(\mathbb{R}^n)\times L^2(\mathbb{R}^n)$, such that, denoting by $w,\tilde{w}$ the solutions of the linear homogeneous Cauchy problems
\begin{align*}
& \partial_t^2 w-\Delta w+b(t) \partial_t w=0, && (w,\partial_t w)(0,\cdot)=(w_0,w_1), \\
& \partial_t^2 \tilde{w}-\Delta \tilde{w}=0, && (\tilde{w},\partial_t \tilde{w})(0,\cdot)=\mathcal{W}(w_0,w_1),
\end{align*} respectively, then
\begin{align*}
\|(w,\partial_t w)(t,\cdot)-(\tilde{w},\partial_t \tilde{w})(t,\cdot)\|_{\mathcal{E}} \to 0 \ \ \mbox{as} \ \ t\to +\infty.
\end{align*}}
 Moreover, in \cite{Wi06} for the linear equation $\partial_t^2 u-\Delta u+b(t)\partial_t u=0$ the same $L^p-L^{p'}$ estimates are prove as for the free wave equation. In a nutshell, we can think to the damping term $b(t)\partial_t u$ as a negligible perturbation when the coefficient $b$ is summable. Concerning the semilinear equation  $\partial_t^2 u-\Delta u+b(t)\partial_t u=|u|^p$, the critical exponent is the same one for the classical semilinear wave equation $\partial_t^2 u-\Delta u=|u|^p$, the so-called Strauss exponent $p_{\mathrm{Str}}(n)$, namely, the positive root of the quadratic equation $(n-1)p^2-(n+1)p-2=0$, cf. \cite{LaiTaka2018,LW202,WY2019}.

Let us review the history and current state of research on the critical curve for weakly coupled system of semilinear wave and damped wave equations.

For the classical weakly coupled system of semilinear wave equations
\begin{equation}\label{weak1}
  \left\{
\begin{aligned}
&\partial_t^2u-\Delta u=|v|^p,&\quad& t>0, \  x\in \mathbb{R}^n,\\
&\partial_t^2v-\Delta v=|u|^q,&\quad& t>0,\  x\in \mathbb{R}^n, \\
&(u,\partial_tu,v,\partial_tv)(0,x)= (u_0,u_1,v_0,v_1)(x),               &\quad& x\in\mathbb{R}^n,
\end{aligned}
 \right.
  \end{equation}
the critical curve in the $p-q$ plane is given by
\begin{align}\label{epq}
\Gamma_{\mathrm{W}}(n,p,q)\doteq  \max\left\{\frac{p+2+q^{-1}}{pq-1},\frac{q+2+p^{-1}}{pq-1}\right\} -\frac{n-1}{2},
\end{align} that is, small data solutions to \eqref{weak1} exist globally (in time) when $\Gamma_{\mathrm{W}}(n,p,q)<0$, while local solutions blow up in finite time if $\Gamma_{\mathrm{W}}(n,p,q) \geq 0$ (under suitable sign assumptions for the Cauchy data).  We refer to  \cite{ AgKuTa2000, De97, DeGeMi1997, DeMi1998, GeTaZh2006, IkSoWa2019, Ku05, KuTa03, KuTaWa2012} for further details. We stress that by saying that the system in \eqref{weak1} is weakly coupled we mean that on the right-hand side of the equation for $u$ (resp. $v$) the nonlinear term depends only on $v$ (resp. $u$). The critical curve in \eqref{epq} generalizes the condition that provides the Strauss exponent, which can be rewritten as $\frac{1+p^{-1}}{p-1}=\frac{n-1}{2}$.

For the weakly coupled system of the classical semilinear wave damped equations
\begin{equation}\label{weak2}
  \left\{
\begin{aligned}
&\partial_t^2u-\Delta u+\partial_tu=|v|^p,&\quad& t>0, \  x\in \mathbb{R}^n,\\
&\partial_t^2v-\Delta v+\partial_tv=|u|^q,&\quad& t>0,\  x\in \mathbb{R}^n,  \\
&(u,\partial_tu,v,\partial_tv)(0,x)= (u_0,u_1,v_0,v_1)(x),               &\quad& x\in\mathbb{R}^n,
\end{aligned}
 \right.
  \end{equation}
combining the results scattered in the literature \cite{Nara2009,Ni2012,NishWaka2014,SunWang2007}, we have that the critical curve is given by
\begin{align}\label{fpq}
\Gamma_{\mathrm{DW}}(n,p,q)\doteq\max\left\{\frac{p+1}{pq-1},\frac{q+1}{pq-1}\right\}-\frac{n}{2}.
\end{align}
This results is a natural generalization of the critical exponent for the single semilinear damped wave equation $\partial_t^2u-\Delta u+\partial_tu=|u|^p$. Indeed,  the critical exponent for this equation (cf. \cite{IT2005,M1976,ToYo2001,Zh2001}) is the same one as for the corresponding semilinear heat equation $\partial_tu-\Delta u=|u|^p$, the celebrated Fujita exponent $p_{\mathrm{Fuj}}(n)\doteq 1+\frac{2}{n}$ (named after the author of \cite{Fujita1969}), and it is obtained by the equation $\frac{1}{p-1}-\frac{n}{2}=0$. It can be shown that as in the case of the single equations, in which $p_{\mathrm{Fuj}}(n)<p_{\mathrm{Str}}(n)$, also for the weakly coupled systems we have that the blow-up region for \eqref{weak1} is broader than the one for \eqref{weak2} due to the presence of the damping terms.

After comparing the critical curves for \eqref{weak1} and \eqref{weak2}, seen how strongly the presence of damping terms modify the blow-up region in the $p-q$ plane, it is natural to wonder how the presence of a damping term in just one wave equation would modify the critical curve. In other words, given a system of a semilinear wave equation and a semilinear classical damped wave equation, which are weakly coupled through the nonolinear terms, namely
\begin{equation}\label{weak3}
  \left\{
\begin{aligned}
&\partial_t^2u-\Delta u+\partial_tu=|v|^p,&\quad& t>0, \  x\in \mathbb{R}^n,\\
&\partial_t^2v-\Delta v=|u|^q,&\quad& t>0,\  x\in \mathbb{R}^n  \\
&(u,\partial_tu,v,\partial_tv)(0,x)= (u_0,u_1,v_0,v_1)(x),               &\quad& x\in\mathbb{R}^n,
\end{aligned}
 \right.
  \end{equation} we want to understand how the single damping term $\partial_t u$ influence the critical curve. The Cauchy problem in \eqref{weak3} is known as Nakao problem in the recent literature \cite{Chen2022,ChenReissig2021,KK2022,PaTa2023,Waka2017} and it is named after the author of \cite{N16,N18}, who first introduced and studied a semilinear coupled system of a heat equation and a wave equation in the case of bounded domains.

A first blow-up result, proved by using the so-called test function method, was obtained by Wakasugi \cite{Waka2017}, who proved the blow-up in finite time of weak solutions to \eqref{weak3} when $(p,q)$ satisfies
\begin{equation}\label{gpq}
\Gamma_{\mathrm{N},1}(n,p,q)\doteq \max\left\{\frac{\frac{q}{2}+1}{pq-1}+\frac{1}{2},\frac{p+1}{pq-1},\frac{q+1}{pq-1}\right\}-\frac{n}{2} \geq 0.
\end{equation}
Roughly speaking, we can say that the test function method relies strongly on the scaling properties of the partial differential operator and it allows to prove sharp blow-up results whenever the critical exponent is determined by these scaling properties (as in the case of the Fujita exponent). For this reason, the fact that the last two components in $\Gamma_{\mathrm{N},1}(n,p,q)$ are exactly the same as in $\Gamma_{\mathrm{DW}}(n,p,q)$ is a consequence of the employment of this blow-up technique. However, for the semilinear wave equation it is known that by using the test function method (in the formulation by Mitidieri-Pohozaev \cite{MPbook}) it is possible to prove the blow-up of weak solutions only for $1<p\leq \frac{n+1}{n-1}$, see also \cite{MP2000}, without reaching the Strauss exponent (which is not related to the scaling properties of the d'Alembertian operator). The exponent $\frac{n+1}{n-1}$ is also called Kato exponent in the framework of the semilinear wave equation with nonnegative power nonlinearity, after the author of \cite{K1980}. Subsequently to Wakasugi's blow-up result, in \cite{ChenReissig2021} Chen-Reissig partially improved the blow-up range for \eqref{weak3}. In \cite{ChenReissig2021} it is used a blow-up technique, which is typically used to deal with the classical semilinear wave equation to obtain as upper bound for the power of the nonlinear term the Strauss exponent, cf. \cite{G1981,K1980,Sc1985,Si1984,Ta2015,Zh2006}. More precisely, by using an iteration argument, the authors of \cite{ChenReissig2021} proved the blow-up of energy solutions to \eqref{weak3} under suitable sign assumptions for the Cauchy data and provided that
\begin{align*}\label{CRconditionNP}
\Gamma_{\mathrm{N},2}(n,p,q)\doteq \max \left\{\frac{2+p^{-1}}{pq-1},\frac{\frac{q}{2}+1}{pq-1},\frac{p+\frac{1}{2}}{pq-1}-\frac{1}{2}\right\}-\frac{n-1}{2}>0.
\end{align*}
In particular, in order to derive a sequence of lower bound estimates for the space averages of the components of a local solution, it is used the slicing procedure from \cite{AgKuTa2000}, in the formulation provided in \cite{ChenPa2020}, so that an unbounded exponential multiplier in the iteration frame can be dealt with. Afterwards, Kita-Kusaba \cite{KK2022} combined the results from \cite{ChenReissig2021,Waka2017} and proved the blow-up of weak solutions to \eqref{weak3} for $\max\left\{\Gamma_{\mathrm{N},1}(n,p,q),\Gamma_{\mathrm{N},2}(n,p,q)\right\}>0$ by using the modified test function method from \cite{IkSoWa2019}. To the best of our knowledge, no global existence results for \eqref{weak3}, that may confirm the optimality of any of the previous blow-up ranges, has been proved in the literature.

Our goal in the present manuscript is to investigate how the presence of the damping term $b(t)\partial_t v$ in \eqref{eqs} modifies the blow-up range in the $p-q$ plane for \eqref{eqs} in comparison to the one for \eqref{weak3} in the scale-invariant case and in the scattering producing case, respectively. On the one hand, for $b(t)=\frac{\mu}{1+t}$ we are going to derive a blow-up range which is obtained by the one for \eqref{weak3} by a translation in the space dimension of magnitude $\mu$. This effect of the scale-invariant damping terms in semilinear wave equations has been observed not only for a power nonlinearity $|u|^p $ \cite{Da2021,AbLuRe2015,IkSo2018,LaiPaTa2025,Pa2024}, but also for a derivative type nonlinearity $|\partial_t u|^p$ \cite{Ben2022,HH2021,PaTu2021} and for a combined nonlineary $|\partial_t u|^p+|u|^q$ \cite{HH2021,LiHeQiao2024}. On the other hand, for $b\in L^1([0,\infty))$ we are going to derive the same blow-up range obtained in \cite{ChenReissig2021} for \eqref{weak3}. The fact that the presence of scattering producing damping terms does not alter the critical exponents (or the critical curves in the case of weakly coupled systems) has been already remarked in the literature \cite{LaiTaka2018,LW202,PaTa2019,WY2019}.

In our approach, we study the evolution in time of the space averages of the components of a local solution to \eqref{eqs}. For this purpose, we apply a comparison argument for a system of ordinary differential inequalities, whose proof is based on an iteration argument.

The remaining sections of this paper are organized as follows: in Section \ref{Section main results}, we state the main results of the paper;  Section 3 is devoted to determining the iteration frames and some lower bound estimates for the space averages of the components of a local solution; in Section \ref{Section comparison arggument ODI}, we state and prove a comparison argument for a weakly coupled system of ordinary differential inequalities and in Section \ref{Section proof of the thms} we apply this result to prove the main theorems; finally, in Section \ref{Section FR&OP}, we discuss some open problems related to \eqref{eqs}.

\noindent\textbf{Notations.}
$a \lesssim b$  means that there exists a constant
$C>0$ independent of $t,x,\varepsilon$ such that $a\leq Cb$. We use the notation $B_R\doteq \{x\in\mathbb{R}^n: \vert x\vert\leq R\}$ for the ball of radius $R>0$ around the origin. By $y_+\doteq \max\{y,0\}$ and $y_-\doteq \max\{-y,0\}$  we denote the positive and the negative part of $y\in\mathbb{R}$, respectively.

\section{Main results} \label{Section main results}
To begin with, we provide the definition of the weak solution to \eqref{eqs} that we are going to use in this paper. We stress that this class of solutions is the broadest class for which our approach can be used; for the sake of simplicity, we still call them weak solutions, even though we are assuming additional regularity with respect to the time variable in comparison to the usual weak solutions.
\begin{definition}\label{engdef1}
Let $\varepsilon>0$ and $(u_0,u_1,v_0,v_1)\in\big(W^{1,1}_{\mathrm{loc}}(\mathbb{R}^n)\times L^1_{\mathrm{loc}}(\mathbb{R}^n)\big)^2$. Let us assume that $u_0,u_1,v_0,v_1$ are compactly supported with 
\begin{equation}\label{initialass1}
\mathrm{supp}(u_0,u_1,v_0,v_1)\subset B_R 
\end{equation} for some $R>0$. We say that $(u,v)$ is a weak solution to the weakly system \eqref{eqs} on $[0,T)$, if
\begin{equation*}
\begin{aligned}
&u\in \mathcal{C}\big([0,T),W^{1,1}_{\mathrm{loc}}(\mathbb{R}^n)\big)\cap \mathcal{C}^1\big([0,T),L^1_{\mathrm{loc}}(\mathbb{R}^n)\big)\cap L^q_{\mathrm{loc}}([0,T)\times\mathbb{R}^n)\\
&v\in \mathcal{C}\big([0,T),W^{1,1}_{\mathrm{loc}}(\mathbb{R}^n)\big)\cap \mathcal{C}^1\big([0,T),L^1_{\mathrm{loc}}(\mathbb{R}^n)\big)\cap L^p_{\mathrm{loc}}([0,T)\times\mathbb{R}^n)
\end{aligned}
\end{equation*}
satisfy the support condition
\begin{equation}\label{support condition (u,v)}
\mathrm{supp}(u,v)(t,\cdot)\subset B_{R+t} \quad \mbox{for any} \ t\in (0,T),
\end{equation}
the initial conditions $u(0,\cdot) = \varepsilon u_0$, $v(0,\cdot) = \varepsilon v_0$ in $L_{\mathrm{loc}}^1(\mathbb{R}^n)$, and the integral equalities
\begin{align}
&\int_0^t\int_{\mathbb{R}^n}(-\partial_su(s,x)\partial_s\phi(s,x)+\partial_su(s,x)\phi(s,x)+\nabla v(s,x)\cdot\nabla\phi(s,x))\mathrm{d}x\, \mathrm{d}s \notag \\
&+\int_{\mathbb{R}^n}\partial_tu(t,x)\phi(t,x)\mathrm{d}x=\varepsilon\int_{\mathbb{R}^n}u_1(x)\phi(0,x)\mathrm{d}x
+\int_0^t\int_{\mathbb{R}^n}\vert v(s,x)\vert^p\phi(s,x)\mathrm{d}x\, \mathrm{d}s \label{inteeq11}
\end{align}
and
\begin{align}
&\int_0^t\int_{\mathbb{R}^n}\big(-\partial_sv(s,x)\partial_s\psi(s,x)+\nabla v(s,x)\cdot\nabla\psi(s,x)+b(s)\partial_sv(s,x)\psi(s,x)\big)\mathrm{d}x\, \mathrm{d}s \notag \\
&+\int_{\mathbb{R}^n}\partial_tv(t,x)\psi(t,x)\mathrm{d}x=\varepsilon\int_{\mathbb{R}^n}v_1(x)\psi(0,x)\mathrm{d}x
+\int_0^t\int_{\mathbb{R}^n}\vert u(s,x)\vert^q\psi(s,x)\mathrm{d}x\, \mathrm{d}s \label{inteeq12}
\end{align}
for any $\phi,\psi\in \mathcal{C}_0^\infty([0,T)\times\mathbb{R}^n)$ and any $t\in(0,T)$.
\end{definition}

\begin{remark}
After further integrations by parts in \eqref{inteeq11} and \eqref{inteeq12}, we get
\begin{align}
&\int_0^t\int_{\mathbb{R}^n}u(s,x)\big(\partial_s^2\phi(s,x)-\Delta\phi(s,x)-\partial_s\phi(s,x)\big)\mathrm{d}x\, \mathrm{d}s \notag \\
&\quad +\int_{\mathbb{R}^n}\big(\partial_tu(t,x)\phi(t,x)-u(t,x)\partial_t\phi(t,x)+u(t,x)\phi(t,x)\big)\mathrm{d}x \notag\\
&=\varepsilon\int_{\mathbb{R}^n}\big((u_0(x)+u_1(x))\phi(0,x)-u_0(x)\partial_t\phi(0,x)\big)\mathrm{d}x+\int_0^t\int_{\mathbb{R}^n}\vert v(s,x)\vert^p\phi(s,x)\mathrm{d}x\, \mathrm{d}s \label{inteeq21}
\end{align}
and
\begin{align}
&\int_0^t\int_{\mathbb{R}^n}v(s,x)\big(\partial_s^2\psi(s,x)-\Delta\psi(s,x)-\partial_s(b(s)\psi(s,x))\big)\mathrm{d}x\, \mathrm{d}s \notag\\
& \quad+\int_{\mathbb{R}^n}\big(\partial_tv(t,x)\psi(t,x)-v(t,x)\partial_t\psi(t,x)+b(t)v(t,x)\psi(t,x)\big)\mathrm{d}x \notag\\
&=\varepsilon\int_{\mathbb{R}^n}\big((b(0) v_0(x)+v_1(x))\psi(0,x)-v_0(x)\partial_t\psi(0,x)\big)\mathrm{d}x+\int_0^t\int_{\mathbb{R}^n}\vert u(s,x)\vert^q\psi(s,x)\mathrm{d}x\, \mathrm{d}s \label{inteeq22}
\end{align}
for any $\phi,\psi\in \mathcal{C}_0^\infty([0,T)\times\mathbb{R}^n)$ and any $t\in(0,T)$
\end{remark}

Let us state our main blow-up results for \eqref{eqs}.

\begin{theorem} \label{theorem1}
Let $b(t)=\frac{\mu}{1+t}$, with $\mu\geq0$ in \eqref{eqs}. Let $(u_0,u_1,v_0,v_1)\in \big(W^{1,1}_{\mathrm{loc}}(\mathbb{R}^n)\times L^1_{\mathrm{loc}}(\mathbb{R}^n)\big)^2$ be nonnegative and nontrivial, and satisfy \eqref{initialass1}.

Let $p,q>1$ such that
\begin{equation}\label{loweralphanmu}
\Gamma(n,p,q,\mu) \doteq \max \left\{\frac{2+p^{-1}}{pq-1}-\frac{n+\mu-1}{2},\frac{q+2}{pq-1}-n+[1-\mu]_+, \frac{2p+1}{pq-1}-n\right\}>0.
\end{equation}
 There exists a $\varepsilon_0=\varepsilon_0(n,\mu,p,q,R,u_0,u_1,v_0,v_1)>0$ such that for any $\varepsilon\in(0,\varepsilon_0]$,  the local (in time) weak solution $(u,v)$ to \eqref{eqs} on $[0,T(\varepsilon))$ blows up in finite time and 
\begin{equation} \label{ub lifespan thm s.i.}
T(\varepsilon)\leq C \varepsilon^{-\tfrac{1}{\Gamma(n,p,q,\mu)}},
\end{equation} where $C>0$ is independent of $\varepsilon$.
\end{theorem} 

\begin{theorem} \label{theorem2}
Let $b\in L^1([0,\infty))$ in \eqref{eqs} be a nonnegative function. Assume that the  data $(u_0,u_1,v_0,v_1)\in \big(W^{1,1}_{\mathrm{loc}}(\mathbb{R}^n)\times L^1_{\mathrm{loc}}(\mathbb{R}^n)\big)^2$ are nonnegative and nontrivial, and satisfy \eqref{initialass1}.
Let $p,q>1$ such that 
\begin{align}\label{blowup range thm s.p.}
\Gamma(n,p,q,0) = \max \left\{\frac{2+p^{-1}}{pq-1}-\frac{n-1}{2}, \frac{q+2}{pq-1}-n+1, \frac{2p+1}{pq-1}-n\right\}>0.
\end{align}

Then there exists a $\varepsilon_0=\varepsilon_0(n,b,p,q,R,u_0,u_1,v_0,v_1)>0$ such that for any $\varepsilon\in(0,\varepsilon_0]$ the local (in time) weak solution $(u,v)$ to the weakly coupled system \eqref{eqs} on $[0,T(\varepsilon))$ blows up in finite time and 
\begin{align} \label{ub lifespan thm s.p.}
T(\varepsilon)\leq C \varepsilon^{-\tfrac{1}{\Gamma(n,p,q,0)}},
\end{align} where $C>0$ is independent of $\varepsilon$.
\end{theorem}

\begin{remark}\label{rem21}
When $\mu\to 0$,  we observe that Theorem \ref{theorem1} reduces to the result for the classical Nakao problem obtained in \cite{ChenReissig2021}.
\end{remark}

\begin{remark}
Theorem \ref{theorem2} shows that when $b\in L^1([0,\infty))$, the blow-up region coincides with that in \cite{ChenReissig2021} for \eqref{weak3}, indicating that the damping term $b(t)\partial_tv$  does not effect the blow-up region in the scattering producing case.
\end{remark}

\section{Iteration frames and first lower bound estimates} \label{sce3}
Let us fix a weak local solution $(u,v)$ to \eqref{eqs} on $[0,T)$, according to Definition \ref{engdef1}. We define the two functionals
\begin{equation*}
U(t)\doteq \int_{\mathbb{R}^n}u(t,x)\, \mathrm{d}x, \qquad  V(t)\doteq \int_{\mathbb{R}^n}v(t,x)\, \mathrm{d}x.
\end{equation*}
Our goal in this section is to derive the iterative framework and lower bound estimates for  $U,V$ in the two cases $b(t)=\frac{\mu}{1+t}$ and $b\in L^1([0,\infty))$, respectively.

Fixed $t\in (0,T)$, choosing as test function $\phi,\psi\in \mathcal{C}_0^\infty\big([0,T)\times\mathbb{R}^n\big)$ such that $\phi,\psi\equiv 1$ on the set $\{(s,x)\in[0,t]\times\mathbb{R}^n:\vert x\vert\leq R+s\}$, from \eqref{inteeq11} and \eqref{inteeq12} we obtain
\begin{equation*}
\begin{aligned}
&\int_0^t\int_{\mathbb{R}^n}\partial_su(s,x)\mathrm{d}x\, \mathrm{d}s+\int_{\mathbb{R}^n}\partial_tu(t,x)\mathrm{d}x=\varepsilon\int_{\mathbb{R}^n}u_1(x)\mathrm{d}x+\int_0^t\int_{\mathbb{R}^n}\vert v(s,x)\vert^p\mathrm{d}x\, \mathrm{d}s\\
\end{aligned}
\end{equation*}
and
\begin{equation*}
\begin{aligned}
&\int_0^t\int_{\mathbb{R}^n}b(s)\partial_sv(s,x)\mathrm{d}x\, \mathrm{d}s+\int_{\mathbb{R}^n}\partial_tv(t,x)\mathrm{d}x=\varepsilon\int_{\mathbb{R}^n}v_1(x)\mathrm{d}x+\int_0^t\int_{\mathbb{R}^n}\vert u(s,x)\vert^q\mathrm{d}x\, \mathrm{d}s.
\end{aligned}
\end{equation*}
By using the functionals $U,V$, the previous two equations can be written as
\begin{align}
& U'(t)+U(t) =U'(0)+U(0)+\int_{0}^t\int_{\mathbb{R}^n}\vert v(s,x)\vert^p\mathrm{d}x\, \mathrm{d}s,\label{ueq}\\
& V'(t)+\int_0^tb(s)V'(s)\mathrm{d}s =V'(0)+\int_0^t\int_{\mathbb{R}^n}\vert u(s,x)\vert^q\mathrm{d}x\, \mathrm{d}s.\label{veq}
\end{align}
Hereafter, we consider separately the case scale-invariant and the case with scattering producing damping term.

\subsection{Scale-invariant damping} \label{itfr}

In this subsection we consider always $b(t)=\frac{\mu}{1+t}$. Multiplying by $\mathrm{e}^t$ both sides of \eqref{ueq}, it results
\begin{equation*}
\frac{\mathrm{d}}{\mathrm{d}t}\big(\mathrm{e}^tU(t)\big)=\big(U'(0)+U(0)\big)\mathrm{e}^t+\mathrm{e}^t\int_{0}^t\int_{\mathbb{R}^n}\vert v(s,x)\vert^p\mathrm{d}x\, \mathrm{d}s.
\end{equation*}
After integrating the previous equation over $[0,t]$, we get
\begin{equation}\label{solutionU}
U(t)=U(0)+U'(0)(1-\mathrm{e}^{-t})+\int_{0}^t\mathrm{e}^{\tau-t}\int_0^\tau\int_{\mathbb{R}^n}\vert v(s,x)\vert^p\mathrm{d}x\, \mathrm{d}s \,\mathrm{d}\tau.
\end{equation}
By differentiating \eqref{veq} with respect to $t$, we find
\begin{equation}\label{diffV}
V^{\prime\prime}(t)+\frac{\mu}{1+t}V^\prime(t)=\int_{\mathbb{R}^n}\vert u(t,x)\vert^q\mathrm{d}x,
\end{equation}
multiplying $(1+t)^\mu$ both sides of  \eqref{diffV}  we obtain
\begin{equation*}
\frac{\mathrm{d}}{\mathrm{d}t}\big((1+t)^\mu V'(t)\big)=(1+t)^\mu\int_{\mathbb{R}^n}\vert u(t,x)\vert^q\mathrm{d}x.
\end{equation*}
Then integrating on $[0,t]$, we arrive at the following representation for $V(t)$:
\begin{align}
V(t)=V(0)&+V'(0)\int_0^t\frac{1}{(1+s)^\mu}\mathrm{d}s 
+ \int_0^t(1+\tau)^{-\mu}\int_0^\tau(1+s)^\mu\int_{\mathbb{R}^n}\vert u(s,x)\vert^q\mathrm{d}x\, \mathrm{d}s\, \mathrm{d}\tau. \label{solutionV}
\end{align}
 By using the support condition \eqref{support condition (u,v)} and H\"older's inequality, we have
\begin{align}
(U(s))^q & \lesssim(R+s)^{n(q-1)}\int_{\mathbb{R}^n}\vert u(s,x)\vert^q \mathrm{d}x,\label{holder1}\\
(V(s))^p & \lesssim(R+s)^{n(p-1)}\int_{\mathbb{R}^n}\vert v(s,x)\vert^q\mathrm{d}x.\label{holder2}
\end{align}
Thanks to the nonnegativity of the initial data, it holds $U(0),U'(0),V(0),V'(0)\geq 0$. Then, plugging the inequality coming from \eqref{holder2} into \eqref{solutionU}, as well as the one from \eqref{holder1} into \eqref{solutionV},  we get the iteration frame
 \begin{align}
&U(t)\geq C\int_0^t\mathrm{e}^{\tau-t}\int_0^\tau(R+s)^{-n(p-1)}(V(s))^p \mathrm{d}s\, \mathrm{d}\tau,\label{U1}\\
&V(t)\geq K\int_0^t(R+\tau)^{-\mu}\int_0^\tau s^\mu(R+s)^{-n(q-1)}(U(s))^q \mathrm{d}s\, \mathrm{d}\tau,\label{V1}
\end{align}
where $C,K$ are positive constants depending only on $n,p,q,R$.

\subsubsection{The first kind of lower bounds for the functionals} \label{subsub scale invariant}
We consider the adjoint equation to the homogeneous equation associated to the $v$-equation in \eqref{eqs}, namely,
\begin{equation}\label{conjugate}
\partial_s^2\Psi(s,x)-\Delta\Psi(s,x)-\partial_s\left(\frac{\mu}{1+s}\Psi(s,x)\right)=0.
\end{equation}
Our aim is to find a solution to \eqref{conjugate}, which is nonnegative and with separated variables.
As $x$-dependent factor we take
\begin{equation}\label{varphi}
 \varphi(x)\doteq 
 \begin{cases}
 \mathrm{e}^{x}+\mathrm{e}^{-x} & \mbox{if} \ n=1, \\
  \displaystyle{\vphantom{\Big(}\int_{\mathbb{S}^{n-1}}\mathrm{e}^{x\cdot \omega}d \sigma_\omega} & \mbox{if} \ n\geqslant 2.
 \end{cases}
\end{equation}
The function $\varphi$  was introduced in \cite{Zh2006} and has been widely used ever since to study the blow-up phenomena for the semilinear wave models. The function $\varphi\in\mathcal{C}^{\infty}(\mathbb{R}^n)$ is positive and satisfies $\Delta \varphi=\varphi$.

 Furthermore, we have the following results, cf. \cite[Equation (2.5)]{Zh2006}.

\begin{lemma}\label{lemma3.3}
Let $r>1$ and $R>0$, then,  there exists $C_0=C_0(n,r,R)>0$ such that
\begin{equation}
\int_{ B_{R+t}}(\varphi(x))^{r}\mathrm{d}x\leq C_0\mathrm{e}^{rt}(R+t)^{(n-1)\left(1-\frac{r}{2}\right)}
\end{equation} for any $t\geq 0$.
\end{lemma}

If we set $\Psi(s,x)=\lambda(s)\varphi(x)$, then $\Psi$ is a solution to \eqref{conjugate} if and only if $\lambda$ solves
\begin{equation}\label{ode}
\frac{\mathrm{d}^2}{\mathrm{d}t^2}\lambda(t)-\frac{\mathrm{d}}{\mathrm{d}t}\left(\frac{\mu}{1+t}\lambda(t)\right)-\lambda(t)=0.
\end{equation}
Let $\mathrm{K}_\ell$ be the modified Bessel function of the second kind of order $\ell\in\mathbb{R}$. Then, $\mathrm{K}_\ell$ satisfies the ODE
\begin{equation}\label{bessel}
z^2\frac{\mathrm{d}^2}{\mathrm{d}z^2}\mathrm{K}_\ell(z)+z\frac{\mathrm{d}}{\mathrm{d}z}\mathrm{K}_\ell(z)-(z^2+\ell^2)\mathrm{K}_\ell(z)=0
\end{equation} and admit the following integral representation formula for $z>0$ (see \cite[Equations (10.25.1) and (10.32.9)]{Tricomi1953})
\begin{equation}\label{defK}
\mathrm{K}_\ell(z)=\int_0^\infty \mathrm{e}^{-z \cosh y}\cosh(\ell y) \mathrm{d}y.
\end{equation}
Moreover, $\mathrm{K}_\ell$ satisfies the following properties
\begin{align}
 \mathrm{K}_\ell(z) &=\sqrt{\frac{\pi}{2z}}\mathrm{e}^{-z}\bigl(1+O(z^{-1})\bigr) \quad \mbox{as} \ z\rightarrow +\infty,\label{prop1}\\
 \frac{\mathrm{d}}{\mathrm{d}z}\mathrm{K}_\ell(z) &=-\mathrm{K}_{\ell+1}(z)+\frac{\ell}{z}\, \mathrm{K}_\ell(z)\label{prop2},
\end{align} cf. \cite[Equations (10.40.2) and (10.29.2)]{Tricomi1953}.
\begin{lemma}\label{consolution}
The function 
\begin{equation} \label{def function lambda}
\lambda(t;\mu)\doteq (1+t)^{\frac{\mu+1}{2}}\mathrm{K}_\ell(1+t)
\end{equation} is a solution to \eqref{ode} for $\ell=\frac{|\mu-1|}{2}$.
\end{lemma}
\begin{proof}
Note that the equation \eqref{ode} is  equivalent to
\begin{equation}\label{equationequal}
(1+t)^2\frac{\mathrm{d}^2\lambda}{\mathrm{d}t^2}(t)-\mu(1+t)\frac{\mathrm{d}\lambda}{\mathrm{d}t}(t)+\big(\mu-(1+t)^2\big)\lambda(t)=0.
\end{equation}
Our goal is to transform \eqref{equationequal} into a modified Bessel equation. Set $\lambda(t)\doteq (t+1)^k \eta(t)$, where $k$ is a positive parameter to be determined later. After some elementary computations, we find that $\lambda$ solves \eqref{equationequal} if and only if $\eta$ solves
\begin{equation*}
(1+t)^2\frac{\mathrm{d}^2 \eta}{\mathrm{d}t^2}(t)+(2k-\mu)(1+t)\frac{\mathrm{d}\eta}{\mathrm{d}t}(t)+\big((k-\mu)(k-1)-(1+t)^2\big)\eta(t)=0.
\end{equation*}
Hence, choosing $k=\frac{\mu+1}{2}$, the previous equation becomes the following modified Bessel equation
\begin{equation*}
(1+t)^2\frac{\mathrm{d}^2 \eta}{\mathrm{d}t^2}(t)+(1+t)\frac{\mathrm{d}\eta}{\mathrm{d}t}(t)-\left(\left(\tfrac{\mu-1}{2}\right)^2+(1+t)^2\right)\eta(t)=0.
\end{equation*}
 Therefore, the function $\lambda$ defined in \eqref{def function lambda} is a solution to \eqref{conjugate} provided that $\ell=\frac{|\mu-1|}{2}$.
\end{proof}

In the next result, we determine the asymptotic behavior of $\lambda$.
\begin{lemma}\label{lambdaprop}
There exists $C_1=C_1(\mu)>0$ such that for any $t\geq0$,
\begin{equation}\label{equalprop}
\frac{1}{C_1}(1+t)^{\frac{\mu}{2}}\mathrm{e}^{-t}\leq\lambda(t;\mu)\leq C_1(1+t)^{\frac{\mu}{2}}\mathrm{e}^{-t}.
\end{equation}
\end{lemma}

\begin{proof}
From \eqref{prop1}, it follows that there exists $T_0>0$ such that
\begin{equation*}
(1+t)^{-\frac{1}{2}}\mathrm{e}^{-t} \lesssim \mathrm{K}_\ell(1+t)\lesssim (1+t)^{-\frac{1}{2}}\mathrm{e}^{-t}
\end{equation*} for any $t>T_0$. Consequently, for any $t>T_0$
\begin{align*}
(1+t)^{\frac{\mu}{2}}\mathrm{e}^{-t} \lesssim \lambda(t;\mu) \lesssim (1+t)^{\frac{\mu}{2}}\mathrm{e}^{-t}.
\end{align*} Since the functions $\lambda(t;\mu)$ and $(1+t)^{-\frac{\mu}{2}}\mathrm{e}^{-t}$ are positive for any $t>0$, it is clear that the previous estimates are satisifed also for $t\in[0,T_0]$ for suitable multiplicative constants.
\end{proof}

Hereafter, we set $\Psi(t,x)\doteq \lambda(t,\mu) \varphi(x)$ with $\lambda(\cdot;\mu)$ given by \eqref{def function lambda}.
We define the auxiliary functional
\begin{equation}\label{def functional V0}
V_0(t)\doteq \int_{\mathbb{R}^n}v(t,x)\Psi(t,x)\, \mathrm{d}x.
\end{equation}
Our next goal is to  derive a lower bound estimate for $U(t)$ with the help of a lower bound estimate for $V_0(t)$.

\begin{lemma}\label{V0es}
Let $v_0,v_1$ be functions satisfying the same assumptions as in the statement of Theorem \ref{theorem1}. Let $V_0$ be the functional associated with a local solution $(u,v)$ of \eqref{eqs} for $b(t)=\frac{\mu}{1+t}$ given by \eqref{def functional V0}.
Then there exists  $D=D(n,\mu,v_0,v_1)>0$ such that 
\begin{align*}
V_0(t)\geq D\varepsilon \qquad \mbox{for any} \ t\in [0,T).
\end{align*} 
\end{lemma}
\begin{proof}
Thanks to the support condition \eqref{support condition (u,v)} for $(u,v)$, we may consider as test function in \eqref{inteeq22} $\Psi\in\mathcal{C}^\infty([0,T)\times \mathbb{R}^n)$. Since $\Psi$ is a solution to the adjoint equation \eqref{conjugate} we get
\begin{align}
&V'_0(t)+\left(\frac{\mu}{1+t}-\frac{2\lambda'(t)}{\lambda(t)}\right)V_0(t)=\varepsilon I[v_0,v_1,\mu]+\int_0^t\int_{\mathbb{R}^n}\vert u(s,x)\vert^q \Psi(s,x) \mathrm{d}x\, \mathrm{d}s, \label{keyinte}
\end{align}
where $$I[v_0,v_1,\mu]\doteq \int_{\mathbb{R}^n}\Big(\lambda(0;\mu)\big(\mu v_0(x)+v_1(x)\big)-\lambda'(0;\mu)v_0(x)\Big)\varphi(x)\mathrm{d}x.$$
Combining \eqref{def function lambda} and \eqref{prop2}, we find that
\begin{align*}
&\frac{\mathrm{d}\lambda}{\mathrm{d}t}(t)=\left(\tfrac{\mu+1}{2}+\ell\right)(1+t)^{\frac{\mu-1}{2}}\mathrm{K}_{\ell}(1+t)-(1+t)^{\frac{\mu+1}{2}}
\mathrm{K}_{\ell+1}(1+t),
\end{align*} where $\ell= \frac{|\mu-1|}{2}$.
Since the modified Bessel functions $\mathrm{K}_{\ell},\mathrm{K}_{\ell+1}$ and the data $v_0,v_1$ are nonnegative functions, we have
\begin{align*}
I[v_0,v_1,\mu]& =\int_{\mathbb{R}^n}\left(\left(\left(\tfrac{\mu-1}{2}-\ell\right) v_0(x)+v_1(x)\right) \mathrm{K}_{\ell}(1)
+\mathrm{K}_{\ell+1}(1)v_0(x)\right)\varphi(x)\mathrm{d}x \\
& \geq \int_{\mathbb{R}^n}\left((1-[1-\mu]_+)v_0(x)+v_1(x)\right)\mathrm{K}_{\ell}(1)\varphi(x)\mathrm{d}x>0,
\end{align*} where we used the monotonicity of $\mathrm{K}_{\ell}$ in $\ell$, indeed, from \eqref{defK} it is clear that $\mathrm{K}_{\ell}(z)$ is increasing with respect to $\ell\geq 0$ for a fixed $z>0$.

Multiplying by $\frac{(1+t)^\mu}{\lambda^2(t)}$  both sides of \eqref{keyinte}, we get
\begin{equation}\label{intekey}
\frac{\mathrm{d}}{\mathrm{d}t}\left(\frac{(1+t)^\mu}{\lambda^2(t)}V_0(t)\right)\geq \varepsilon I[v_0,v_1,\mu]  \,\frac{(1+t)^\mu}{\lambda^2(t)}.
\end{equation}
Integrating  \eqref{intekey} over the interval $[0, t]$, we obtain that  for any $t
\geq 0$,
\begin{equation*}
\begin{aligned}
V_0(t) & \geq \frac{V_0(0)}{\lambda^2(0)}\frac{\lambda^2(t)}{(1+t)^\mu}+\varepsilon I[v_0,v_1,\mu]\frac{\lambda^2(t)}{(1+t)^\mu}\int_0^t\frac{(1+s)^\mu}{\lambda^2(s)}\mathrm{d}s\\
&\geq \frac{V_0(0)}{C_1^{4}}  \mathrm{e}^{-2t} + \frac{\varepsilon I[v_0,v_1,\mu]}{C_1^{4}}\mathrm{e}^{-2t}\int_{0}^t\mathrm{e}^{2s}\mathrm{d}s\\
& = \frac{\varepsilon}{C_1^{4}} \left(\lambda(0;\mu) \int_{\mathbb{R}^n} v_0(x) \varphi(x)\mathrm{d}x \, \mathrm{e}^{-2t} + \frac{1}{2} I[v_0,v_1,\mu] \left(1-\mathrm{e}^{-2t}\right) \right) \geq D\varepsilon,
\end{aligned}
\end{equation*}
where $D\doteq \frac{1}{2C_1^4}\min\left\{2\lambda(0;\mu) \|v_0 \, \varphi\|_{L^1(\mathbb{R}^n)},  I[v_0,v_1,\mu] \right\}$.
\end{proof}

Next, we derive the lower bound estimate for $U$.
By Lemma  \ref{lemma3.3} and Lemma  \ref{lambdaprop}  we have
\begin{equation*}
\begin{aligned}
\int_{B_{R+t}}\left(\Psi(t,x)\right)^{\frac{p}{p-1}}\mathrm{d}x&\leq C_0C_1^{\frac{p}{p-1}}(1+t)^{\frac{\mu p}{2(p-1)}}\mathrm{e}^{-\frac{p}{p-1}t}\mathrm{e}^{\frac{p}{p-1}t}(R+t)^{(n-1)\big(1-\frac{p}{2(p-1)}\big)}\\
&\leq C_2 (R+t)^{\frac{\mu p}{2(p-1)}+(n-1)\big(1-\frac{p}{2(p-1)}\big)},
\end{aligned}
\end{equation*}
for some $C_2=C_2(n,\mu,p,R)>0$.
Then, by H\"older's inequality and \eqref{support condition (u,v)},
\begin{align}
\int_{\mathbb{R}^n}\vert v(t,x)\vert^p \mathrm{d}x&\geq (V_0(t))^p\left(\int_{B_{R+t}}\left(\Psi(t,x)\right)^{\frac{p}{p-1}}\mathrm{d}x\right)^{-(p-1)}\notag \\
&\geq C_3\varepsilon^p(R+t)^{n-1-\frac{n+\mu-1}{2}p}, \label{holder3}
\end{align}
where $C_3\doteq D^pC_2^{-(p-1)}$. Combining \eqref{holder3} and \eqref{solutionU},  we obtain
\begin{align}
U(t) & \geq C_3\varepsilon^{p}\mathrm{e}^{-t}\int_0^t\mathrm{e}^\tau\int_0^\tau(R+s)^{n-1-\frac{n+\mu-1}{2}p}\mathrm{d}s\, \mathrm{d}\tau \notag \\
&\geq C_3\varepsilon^{p}\mathrm{e}^{-t}(R+t)^{-\frac{n+\mu-1}{2}p}\int_0^t\mathrm{e}^\tau\int_0^\tau s^{n-1}\mathrm{d}s\, \mathrm{d}\tau \notag\\
&\geq \frac{C_3\varepsilon^{p}}{n}(R+t)^{-\frac{n+\mu-1}{2}p}\, \mathrm{e}^{-t}\int_{\frac{t}{2}}^t\mathrm{e}^\tau\tau^n \mathrm{d}\tau \notag \\
 &\geq \frac{C_3\varepsilon^{p}}{n2^n}(R+t)^{-\frac{n+\mu-1}{2}p}\, t^n\mathrm{e}^{-t} \int_{\frac{t}{2}}^t\mathrm{e}^\tau \mathrm{d}\tau \notag \\
&\geq  C_4 \varepsilon^{p}(R+t)^{-\frac{n+\mu-1}{2}p}\, t^n \label{1st lb U}
\end{align} for $t\in[1,T)$, where $C_4\doteq \frac{C_3}{n2^n} (1-\mathrm{e}^{-\frac{1}{2}})$.

\subsubsection{The second kind of lower bounds for the functionals}
From the representation  of $V(t)$ in \eqref{solutionV}, we have
\begin{align}
V(t)&\geq V'(0)\int_0^t\frac{\mathrm{d}s}{(1+s)^\mu} 
\geq C_5 \varepsilon 
\begin{cases}
1, & \mbox{if} \ \mu >1,\\
\ln(1+t) & \mbox{if} \ \mu =1,\\
t^{1-\mu} &  \mbox{if} \ \mu \in [0,1).
\end{cases} \label{lb V}
\end{align}
for $t\in [1,T)$, where $C_5=C_5(\mu,v_1)>0$. Analogously, from the representation of $U(t)$ in \eqref{solutionU}, it follows immediately that
\begin{align} \label{2nd lb U}
U(t)\geq U(0)= C_6 \varepsilon
\end{align} for $t\in [0,T)$, where $C_6\doteq \|u_0\|_{L^1(\mathbb{R}^n)}$.

\subsection{Scattering producing damping}\label{loeboun}
 In this case $b\in L^1([0,\infty))$. In order to derive the iteration frame for $V(t)$, inspired by \cite{LaiTaka2018}, we introduce the function
\begin{equation}\label{defm}
m(t)\doteq \mathrm{e}^{-\int_t^\infty b(s)\mathrm{d}s}.
\end{equation}
Since $b$ is nonnegative and summable, it follows that $m$ is increasing and bounded, so
\begin{equation}\label{propm}
0<\mathrm{e}^{-\|b\|_{L^1(\mathbb{R}_+)}}=m(0)\leq m(t)\leq1, 
\end{equation}
for any $t\geqslant 0$ and, moreover,
\begin{align} \label{propm 2}
 m^\prime(t)=b(t)m(t) .
\end{align}
By differentiating \eqref{veq}  with respect to $t$, we have
\begin{equation}\label{diffv}
V^{\prime\prime}(t)+b(t)V^\prime(t)=\int_{\mathbb{R}^n}\vert u(t,x)\vert^q \mathrm{d}x.
\end{equation}
Multiplying  both sides of  \eqref{diffv} by $m(t)$, we obtain
\begin{equation}\label{diffv2}
\frac{\mathrm{d}}{\mathrm{d}t}\big(m(t)V^\prime(t)\big)=m(t)\int_{\mathbb{R}^n}\vert u(t,x)\vert^q \mathrm{d}x.
\end{equation}
Integrating  \eqref{diffv2} over the interval $[0,t]$  and using \eqref{propm}, we find
\begin{align}
V^\prime(t)&\geq m(t)V^\prime(t)\geq m(t)V^\prime(t)-m(0)V^\prime(0) \notag \\
&=\int_0^tm(s)\int_{\mathbb{R}^n}\vert u(s,x)\vert^q\mathrm{d}x\, \mathrm{d}s \geq m(0)\int_0^t\int_{\mathbb{R}^n}\vert u(s,x)\vert^q \mathrm{d}x\, \mathrm{d}s. \label{diffv3}
\end{align}
By integrating \eqref{diffv3} over $[0,t]$ and using \eqref{holder1} along with the nonnegativity of the initial data, we obtain the following estimate for $V(t)$:
\begin{equation}\label{framev}
V(t)\geq K\int_0^t\int_0^s(R+\tau)^{-n(q-1)}(U(\tau))^q \mathrm{d}\tau \,\mathrm{d}s,
\end{equation}
where $K=K(n,R,b)>0$.
Clearly, we can repeat the same computations as in Subsection \ref{itfr} to estimate $U(t)$ from below, obtaining again the inequality in \eqref{U1}. Summarizing, in the scattering producing case the iteration frame for $(U,V)$ is given by \eqref{U1} and \eqref{framev}.

\subsubsection{The first kind of lower bounds for the functionals} Let $\Phi(t,x)\doteq \mathrm{e}^{-t}\varphi(x)$, where $\varphi$ is defined in \eqref{varphi}. We introduce the auxiliary functional
\begin{equation}\label{Phi}
V_1(t)\doteq \int_{\mathbb{R}^n}v(t,x)\Phi(t,x)\mathrm{d}x.
\end{equation}
Analogously to what we did in Subsection \ref{subsub scale invariant} with $V_0$, now we are going to determine a lower bound estimate for $V_1$, that will in turn provide  a first lower bound for $U$.

\begin{lemma}\label{propv1}
Let $v_0,v_1$ be functions satisfying the same assumptions as in the statement of Theorem \ref{theorem2}. Let $V_1$ be the functional associated with a local solution $(u,v)$ of \eqref{eqs} for $b\in L^1([0,\infty))$ given by \eqref{Phi}.
There exists a $\tilde{D}=\tilde{D}(n,b,v_0,v_1)>0$ such that,
\begin{equation}\label{lowerv}
V_1(t)\geq\tilde{D}\varepsilon \qquad \mbox{for any} \ t\in [0,T).
\end{equation}
\end{lemma}

\begin{proof}
We consider $\Phi$ as test function into \eqref{inteeq12}. Differentiating the resulting equation, we get
\begin{align}
&\frac{\mathrm{d}}{\mathrm{d}t}\int_{\mathbb{R}^n}\partial_tv(t,x)\Phi(t,x)\mathrm{d}x-\int_{\mathbb{R}^n}v(t,x)\Delta\Phi(t,x)\mathrm{d}x-
\int_{\mathbb{R}^n}\partial_tv(t,x)\partial_t\Phi(t,x)\mathrm{d}x \notag \\
& \quad+ b(t)\int_{\mathbb{R}^n}\partial_tv(t,x)\Phi(t,x)\mathrm{d}x=\int_{\mathbb{R}^n}\vert u(t,x)\vert^q\Phi(t,x)\mathrm{d}x. \label{Phiv}
\end{align}
Multiplying both sides of \eqref{Phiv} by $m(t)$ and using \eqref{propm}, it results
\begin{align}
&\frac{\mathrm{d}}{\mathrm{d}t}\Big(m(t)\int_{\mathbb{R}^n}\partial_tv(t,x)\Phi(t,x)\mathrm{d}x\Big)+m(t)\int_{\mathbb{R}^n}
\partial_tv(t,x)\Phi(t,x)\mathrm{d}x \notag \\
&-m(t)\int_{\mathbb{R}^n}v(t,x)\Phi(t,x)\mathrm{d}x=m(t)\int_{\mathbb{R}^n}\vert u(t,x)\vert^q\Phi(t,x)\mathrm{d}x, \label{Phiv2}
\end{align}
where we used $\partial_t \Phi=-\Phi$ and $\Delta\Phi=\Phi$. Since 
\begin{align*}
V_1'(t) & = \int_{\mathbb{R}^n} \partial_t v(t,x) \Phi(t,x) \mathrm{d}x -  \int_{\mathbb{R}^n} v(t,x) \Phi(t,x) \mathrm{d}x \\
& = \int_{\mathbb{R}^n} \partial_t v(t,x) \Phi(t,x) \mathrm{d}x -  V_1(t),
\end{align*} we may rewrite \eqref{Phiv2} as follows:
\begin{align*}
&\frac{\mathrm{d}}{\mathrm{d}t}\left(m(t)(V'_1(t)+V_1(t))\right)+m(t)V'_1(t)=m(t)\int_{\mathbb{R}^n}\vert u(t,x)\vert^q\Phi(t,x)\mathrm{d}x. 
\end{align*}
Integrating the previous equation over $[0,t]$ and using the definition of $V_1(t)$, we get
\begin{align*}
&m(t)\big(V_1^\prime(t)+V_1(t)\big)-m(0)\varepsilon\int_{\mathbb{R}^n}v_1(x)\varphi(x)\mathrm{d}x+\int_0^tm(s)V_1^\prime(s)\mathrm{d}s =\int_0^tm(s)\int_{\mathbb{R}^n}\vert u(s,x)\vert^q\Phi(s,x)\mathrm{d}x\mathrm{d}s. 
\end{align*}
Integrating by parts and using \eqref{propm 2}, we find
\begin{align}
&m(t)\big(V_1^\prime(t)+2V_1(t)\big)-\varepsilon m(0)J[v_0,v_1] =\int_0^t b(s) m(s)V_1(s)\mathrm{d}s+\int_0^tm(s)\int_{\mathbb{R}^n}\vert u(s,x)\vert^q\Phi(s,x)\mathrm{d}x\mathrm{d}s, \label{meq1}
\end{align}
where $$ J[v_0,v_1]\doteq \int_{\mathbb{R}^n}\big(v_0(x)+v_1(x)\big)\varphi(x)\mathrm{d}x\geq0.$$ Employing \eqref{propm}, from \eqref{meq1} we conclude that
\begin{align}
V_1^\prime(t)+2V_1(t) &\geq \varepsilon \frac{m(0)}{m(t)}J[v_0,v_1]+\frac{1}{m(t)}\int_0^tm(s)b(s)V_1(s)\mathrm{d}s \notag \\
&\geq \varepsilon m(0) J[v_0,v_1]+m(0)\int_0^t b(s)V_1(s)\mathrm{d}s.
\label{meq2}
\end{align}
Since $V_1'(t)+2V_1(t)=\mathrm{e}^{-2t}\frac{\mathrm{d}}{\mathrm{d}t}(\mathrm{e}^{2t}V_1(t))$, multiplying both sides of \eqref{meq2} by $\mathrm{e}^{2t}$ and integrating the resulting inequality, we obtain
\begin{align}\label{meq3}
&\mathrm{e}^{2t}V_1(t)-\varepsilon\int_{\mathbb{R}^n}v_0(x)\varphi(x)\mathrm{d}x\geq \varepsilon \frac{m(0)}{2} J[v_0,v_1](\mathrm{e}^{2t}-1) + m(0)\int_0^t \mathrm{e}^{2s}\int_0^s b(\tau)V_1(\tau)\mathrm{d}\tau \, \mathrm{d}s.
\end{align}
By \eqref{meq3}, it follows that
\begin{align*}
V_1(t) &\geq \varepsilon \mathrm{e}^{-2t} \int_{\mathbb{R}^n}v_0(x)\varphi(x)\mathrm{d}x + \varepsilon \frac{m(0)}{2} J[v_0,v_1](1-\mathrm{e}^{-2t}) + m(0)\mathrm{e}^{-2t}\int_0^t \mathrm{e}^{2s}\int_0^s b(\tau)V_1(\tau)\mathrm{d}\tau \, \mathrm{d}s \\
&\geq  m(0)\mathrm{e}^{-2t}\int_0^t \mathrm{e}^{2s}\int_0^s b(\tau)V_1(\tau)\mathrm{d}\tau \, \mathrm{d}s.
\end{align*} Using an elementary contradiction argument, we see that $V_1$ is a nonnegative function as long as it exists. Consequently, we have
\begin{align}
V_1(t) & \geq \varepsilon \mathrm{e}^{-2t} \int_{\mathbb{R}^n}v_0(x)\varphi(x)\mathrm{d}x + \varepsilon \frac{m(0)}{2} (1-\mathrm{e}^{-2t}) \int_{\mathbb{R}^n}\big(v_0(x)+v_1(x)\big)\varphi(x)\mathrm{d}x \notag\\
& \geq \tilde{D}\varepsilon, \label{esv1t}
\end{align} where $\tilde{D}\doteq \min\left\{ \|v_0\, \varphi\|_{L^1(\mathbb{R}^n)},\frac{m(0)}{2}J[v_0,v_1] \right\}$.
\end{proof}

By H\"older's inequality and Lemmas \ref{lemma3.3}, \ref{propv1}, it follows that
\begin{align}
\int_{\mathbb{R}^n}\vert v(t,x)\vert^p \mathrm{d}x & \geq (V_1(t))^p\Big(\int_{\mathbb{R}^n}\left(\Phi(t,x)\right)^{\frac{p}{p-1}}\mathrm{d}x\Big)^{-(p-1)} \notag \\ &\geq C_7\varepsilon^p(R+t)^{(n-1)\left(1-\frac{p}{2}\right)}, \label{holderv1}
\end{align}
where $C_7\doteq \tilde{D}^pC_0^{-(p-1)}>0$. Combining \eqref{ueq} and \eqref{holderv1}, we have
\begin{align}
U(t) &\geq C_7\varepsilon^p\mathrm{e}^{-t}\int_0^t\mathrm{e}^\tau\int_0^\tau(R+s)^{(n-1)\left(1-\frac{p}{2}\right)}\mathrm{d}s\, \mathrm{d}\tau \notag\\
&\geq C_7\varepsilon^p\mathrm{e}^{-t}(R+t)^{-\frac{(n-1)p}{2}}\int_0^t\mathrm{e}^\tau\int_0^\tau s^{n-1}\mathrm{d}s\, \mathrm{d}\tau \notag \\
&\geq \frac{C_7}{n}\varepsilon^p\mathrm{e}^{-t}(R+t)^{-\frac{(n-1)p}{2}}\int_{\frac{t}{2}}^t\mathrm{e}^\tau\tau^n \mathrm{d}\tau \notag\\
&\geq C_8 \varepsilon^p(R+t)^{-\frac{(n-1)p}{2}}\, t^n. \label{lowerU11}
\end{align} for any $t\in[1,T)$, where $C_8\doteq \frac{C_7}{n2^n} (1-\mathrm{e}^{-\frac{1}{2}})$.

\subsubsection{The second kind of lower bounds for the functionals}
From \eqref{diffv2} it follows that
\begin{align*}
V(t)=V(0)+m(0)V^\prime(0)\int_0^t\frac{ds}{m(s)}+\int_0^t\frac{1}{m(s)}\int_0^sm(\tau)\int_{\mathbb{R}^n}\vert u(\tau,x)\vert^q\mathrm{d}x\, \mathrm{d}\tau \, \mathrm{d}s.
\end{align*}
Then, by \eqref{propm} we get
\begin{align}
V(t)\geq V(0)+ m(0)V^\prime(0)\int_0^t\frac{\mathrm{d}s}{m(s)}\geq V(0)+ V^\prime(0)t \geq C_8\varepsilon (1+t). \label{lowerV11}
\end{align} for any $t\in[0,T)$, where $C_8\doteq \min\{\|v_0\|_{L^1(\mathbb{R}^n)},\|v_1\|_{L^1(\mathbb{R}^n)}\}$.

Finally, as in Subsection \ref{subsub scale invariant}, the lower bound estimate for $U$ in 
\eqref{2nd lb U} is satisfied.

\section{Comparison argument for a weakly coupled system of ODI} \label{Section comparison arggument ODI}

In this section, we state and prove a comparison argument for a weakly coupled system of \emph{ordinary differential inequalities} (ODI). This result is a generalization, in some sense, of the prototype comparison argument for the classical semilinear wave equation in \cite[Lemma 2.1]{Ta2015}. In order to prove this comparison argument, we are going to use the corresponding equivalent integral formulation that allow us to determine a sequence of lower bound estimates for one of the components of a solution to the system of ODI. As the parameter indexing this sequence goes to $\infty$, we see that one of the components cannot be finite for all times and, as byproduct of this procedure, we derive an upper bound estimate for the lifespan for the corresponding component as well.

In the iteration argument, we use a slicing procedure to handle an unbounded exponential multiplier, following the approach introduced in \cite[Section 4]{ChenPa2020}. Moreover, in the proof we combine some ideas from \cite{LaiTaka2018,LaiTaWa2017,PaTa2023}.

\begin{proposition} \label{Prop comparison argument} Let either $b(t)=\dfrac{\mu}{1+t}$ with $\mu>0$ or $b\in L^1([0,\infty))$ be a nonnegative function. Let $F,G\in\mathcal{C}^2([0,T))$ satisfy the following weakly coupled systems of ODI
\begin{align}\label{wcs ODI}
\begin{cases}
F''(t)+F'(t)\geq B (R+t)^{-r} |G(t)|^p & \mbox{for any} \ t\in [0,T), \\
G''(t)+b(t) G'(t)\geq \widetilde{B} (R+t)^{-\rho} |F(t)|^q & \mbox{for any} \ t\in [0,T), 
\end{cases}
\end{align} where $B,\widetilde{B},R>0$, $p,q>1$ and $r,\rho \in\mathbb{R}$ are fixed parameters.

Assume that $F(0),F'(0),G(0),G'(0)\geq 0$ and that at least one between the following lower bound estimates is satisfied
\begin{align}
F(t)\geq A \, t^{a} & \qquad \mbox{for any} \ t\in [T_0,T), \label{lb F prop}\\
G(t)\geq \widetilde{A}\, t^{\alpha} & \qquad \mbox{for any} \ t\in [T_0,T), \label{lb G prop}
\end{align} where $T_0,A,\widetilde{A}>0$ and $a,\alpha \in\mathbb{R}$ are fixed parameters such that
\begin{align}\label{bu range prop}
\max\left\{a+\frac{2p+1-(\rho p+r)}{pq-1},\alpha+\frac{q+2-(rq+\rho)}{pq-1}\right\}>0.
\end{align}

\noindent \textsc{(1)} There exists $A_0=A_0(a,T_0,R,p,q,r,\rho)>0$ such that if \eqref{lb F prop} holds for a multiplicative constant $A\in (0,A_0]$ and $a(pq-1)+2p+1>\rho p+r$, then $F$ blows up in finite time and 
\begin{align}
T\leq C A^{-\big(a+\frac{2p+1-(\rho p+r)}{pq-1}\big)^{-1}}, \label{ub est F}
\end{align} where $C$ is a positive constant independent of $A$. 

\noindent \textsc{(2)} There exists $\widetilde{A}_0=\widetilde{A}_0(\alpha,T_0,R,p,q,r,\rho)>0$ such that if \eqref{lb G prop} holds for a multiplicative constant $\widetilde{A}\in (0,\widetilde{A}_0]$ and $\alpha(pq-1)+q+2>rq+\rho$, then $G$ blows up in finite time and 
\begin{align}
T\leq \widetilde{C} \widetilde{A}^{-\big(\alpha+\frac{q+2-(rq+\rho)}{pq-1}\big)^{-1}}, \label{ub est G}
\end{align} where $\widetilde{C}$ is a positive constant independent of $\widetilde{A}$.
\end{proposition}

\begin{remark} We underline that the smallness assumption for the multiplicative constants $A,\widetilde{A}$ is crucial to get the upper bound estimates for the lifespan in \eqref{ub est F} and in \eqref{ub est G}, respectively. Of course, if the multiplicative constant is not small, the blow-up result is still true. Furthermore, if we prove the blow-up in finite time of $F$ or $G$, \eqref{ub est F} or \eqref{ub est G} are obviously fulfilled whenever $A$ or $\widetilde{A}$ runs in a compact interval away from $0$, upon a modification of the multiplicative constant $C$. 
\end{remark}

\begin{proof}
We are going to prove the results separately in the scale-invariant case $b(t)=\frac{\mu}{1+t}$ and in the scattering producing case $b\in L^1([0,\infty))$.
\paragraph*{\textit{Scale invariant case}}
Let us begin by rewriting \eqref{wcs ODI} in an integral form. Using the identities
\begin{align*}
&F''(t)+F'(t)  = \frac{\mathrm{d}}{\mathrm{d}t}\left(\mathrm{e}^{-t} \frac{\mathrm{d}}{\mathrm{d}t}\left(\mathrm{e}^{t}F(t)\right)\right),\\
&G''(t)+\frac{\mu}{1+t}G'(t) = (1+t)^{-\mu}\frac{\mathrm{d}}{\mathrm{d}t}\left((1+t)^{\mu}G'(t)\right),
\end{align*} the nonegativity of the initial conditions $F(0),F'(0),G(0),G'(0)$ and repeating similar computations to those in in Section \ref{sce3}, we have
\begin{align}
F(t)&\geq B \,\mathrm{e}^{-t}\int_0^t\mathrm{e}^{\tau}\int_0^\tau(R+s)^{-r}|G(s)|^p \mathrm{d}s\, \mathrm{d}\tau,\label{IF F}\\
G(t)&\geq K\int_0^t(R+\tau)^{-\mu}\int_0^\tau s^\mu(R+s)^{-\rho}|F(s)|^q \mathrm{d}s\, \mathrm{d}\tau,\label{IF G}
\end{align} where $K=K(\widetilde{B},R,\mu)>0$. From \eqref{IF F}, \eqref{IF G} it follows that $F,G$ are nonnegative functions. 

\noindent (1) Let us consider the case $a(pq-1)+2p+1>\rho p+r$ provided that \eqref{lb F prop} holds.

The inequalities \eqref{IF F}-\eqref{IF G} will serve as iteration frame to get a sequence of lower bound estimates for $F$ starting from the one in \eqref{lb F prop}.

 Let us define $\ell_0\doteq \max\{1,\frac{1}{T_0}\}$ and $\ell_k\doteq 1+(pq)^{-k}$ for any $k\in\mathbb{N}, k \geq 1$. We denote by $\{L_j\}_{j\in\mathbb{N}}$ the sequence of the partial products of the convergent infinite product $\prod_{k=0}^\infty \ell_k$, namely, $$L_j\doteq \prod_{k=0}^{j} \ell_k$$ for any $j\in\mathbb{N}$. The strictly increasing sequence $\{L_j\}_{j\in\mathbb{N}}$ will be used in the slicing procedure to handle the factor $\mathrm{e}^\tau$ in \eqref{IF F}. Moreover, we denote
\begin{align*}
L_\infty \doteq \prod_{k=0}^\infty \ell_k \equiv \sup_{j\in\mathbb{N}} L_j.
\end{align*}

Let us prove that for any $j\in\mathbb{N}$ the following lower bound estimate for $F$ holds:
\begin{align}
F(t)\geq B_j (R+t)^{-a_j} (t-L_jT_0)^{b_j} \qquad \mbox{for any} \, t\in[L_jT_0,T), \label{lb Fj}
\end{align} where $\{B_j\}_{j\in\mathbb{N}}$, $\{a_j\}_{j\in\mathbb{N}}$, $\{b_j\}_{j\in\mathbb{N}}$ are sequences of nonnegative real numbers to be determined later.

We are going to prove \eqref{lb Fj} by induction on $j\in\mathbb{N}$. For $j=0$ \eqref{lb Fj} is true thanks to \eqref{lb F prop}, provided that $B_0\doteq A$, $a_0\doteq a_-$, $b_0\doteq a_+$. Let us assume that \eqref{lb Fj} is satisfied for some $j\in\mathbb{N}$, then, we have to prove that \eqref{lb Fj} holds for $j+1$ as well. 

Combining \eqref{lb Fj} and \eqref{IF G}, we have for $t\in [L_jT_0,T)$
\begin{align}
G(t) &\geq K\int_{L_jT_0}^t(R+\tau)^{-\mu}\int_{L_jT_0}^\tau s^\mu(R+s)^{-\rho}(F(s))^q \mathrm{d}s\, \mathrm{d}\tau \notag \\
&\geq K (R+t)^{-\mu-\rho_+} \int_{L_jT_0}^t\int_{L_jT_0}^\tau s^{\mu+\rho_-}(F(s))^q \mathrm{d}s\, \mathrm{d}\tau \notag \\
&\geq K B_j^q (R+t)^{-(\mu+\rho_++q a_j)} \int_{L_jT_0}^t\int_{L_jT_0}^\tau (s-L_jT_0)^{\mu+\rho_-+qb_j} \mathrm{d}s\, \mathrm{d}\tau \notag \\
&\geq K B_j^q (\mu+\rho_-+qb_j+1)^{-1}(\mu+\rho_-+qb_j+2)^{-1}(R+t)^{-(\mu+\rho_++q a_j)} (t-L_jT_0)^{\mu+\rho_-+qb_j+2} \notag \\
&\geq K B_j^q (\mu+\rho_-+qb_j+2)^{-2}(R+t)^{-(\mu+\rho_++q a_j)} (t-L_jT_0)^{\mu+\rho_-+qb_j+2}. \label{lb G inter}
\end{align}

Next, plugging the lower bound for $G$ in \eqref{lb G inter} into \eqref{IF F}, for $t\in[L_jT_0,T)$ we get
\begin{align}
F(t) &\geq B \,\mathrm{e}^{-t}\int_{L_jT_0}^t\mathrm{e}^{\tau}\int_{L_jT_0}^\tau(R+s)^{-r}(G(s))^p \mathrm{d}s\, \mathrm{d}\tau \notag \\
&\geq B  (R+t)^{-r_+}\, \mathrm{e}^{-t}\int_{L_jT_0}^t\mathrm{e}^{\tau}\int_{L_jT_0}^\tau(s-L_jT_0)^{r_-}(G(s))^p \mathrm{d}s\, \mathrm{d}\tau \notag \\
&\geq B K^p B_j^{pq} (\mu+\rho_-+qb_j+2)^{-2p} (R+t)^{-[r_++(\mu+\rho_++q a_j)p]}\, \mathrm{e}^{-t} \notag \\ & \qquad \times \int_{L_jT_0}^t\mathrm{e}^{\tau}\int_{L_jT_0}^\tau(s-L_jT_0)^{r_-+(\mu+\rho_-+qb_j+2)p} \mathrm{d}s\, \mathrm{d}\tau \notag \\
&\geq B K^p B_j^{pq} (\mu+\rho_-+qb_j+2)^{-2p} (r_-+(\mu+\rho_-+qb_j+2)p+1)^{-1}  \notag \\ & \qquad \times  (R+t)^{-[r_++(\mu+\rho_++q a_j)p]}\, \mathrm{e}^{-t} \int_{L_jT_0}^t\mathrm{e}^{\tau}(\tau-L_jT_0)^{r_-+(\mu+\rho_-+qb_j+2)p+1} \mathrm{d}\tau. \label{lb F inter}
\end{align}
Let us estimate the last $\tau$-integral for $t\in[L_{j+1}T_0,T)$: since $t\geq L_{j+1}T_0= \ell_{j+1} L_j T_0$ we may shrink the integration interval from $[L_jT_0,t]$ to $\left[\frac{t}{\ell_{j+1}},t\right]$, consequently,
\begin{align}
& \mathrm{e}^{-t} \int_{L_jT_0}^t\mathrm{e}^{\tau}(\tau-L_jT_0)^{r_-+(\mu+\rho_-+qb_j+2)p+1}  \notag \\
& \quad \geq \mathrm{e}^{-t} \int_{\tfrac{t}{\ell_{j+1}}}^t\mathrm{e}^{\tau}(\tau-L_jT_0)^{r_-+(\mu+\rho_-+qb_j+2)p+1}  \mathrm{d}\tau \notag \\
& \quad \geq \ell_{j+1}^{-(r_-+(\mu+\rho_-+qb_j+2)p+1)} (t-\ell_{j+1} L_jT_0)^{r_-+(\mu+\rho_-+qb_j+2)p+1}  \mathrm{e}^{-t} \int_{\tfrac{t}{\ell_{j+1}}}^t\mathrm{e}^{\tau}  \mathrm{d}\tau \notag \\
& \quad = \ell_{j+1}^{-(r_-+(\mu+\rho_-+qb_j+2)p+1)} (t-L_{j+1}T_0)^{r_-+(\mu+\rho_-+qb_j+2)p+1}  \left(1-\mathrm{e}^{-\big(1-\frac{1}{\ell_{j+1}}\big)t}\right) \label{tau integral}.
\end{align}
Notice that for $t\geq L_{j+1}T_0$ it holds
\begin{align}
1-\mathrm{e}^{-\big(1-\frac{1}{\ell_{j+1}}\big)t} & = 1-\mathrm{e}^{-\frac{\ell_{j+1}-1}{\ell_{j+1}}t}\geq 1-\mathrm{e}^{-\frac{\ell_{j+1}-1}{\ell_{j+1}}L_{j+1}T_0} \notag \\
 & = 1-\mathrm{e}^{-(\ell_{j+1}-1)L_{j}T_0} \geq  1-\mathrm{e}^{-(\ell_{j+1}-1)L_{0}T_0} \notag \\
 & = 1-\mathrm{e}^{-(\ell_{j+1}-1)} \geq 1-\left(1-(\ell_{j+1}-1)+\tfrac{1}{2}(\ell_{j+1}-1)^2\right) \notag \\
 & = (\ell_{j+1}-1)\left(1-\tfrac 12 (\ell_{j+1}-1)\right) = (pq)^{-2j} \left((pq)^j-\tfrac 12\right) \notag \\
 & \geq  (pq)^{-2(j+1)} \left(pq-\tfrac 12\right), \label{exp factor}
\end{align} where in the third inequality we used $\mathrm{e}^{-\sigma}\leq 1-\sigma+\frac{\sigma^2}{2}$ for $\sigma\geq 0$.

Combining \eqref{lb F inter}, \eqref{tau integral} and \eqref{exp factor}, we obtain
\begin{align}
F(t) &\geq \frac{B K^p B_j^{pq} \left(pq-\tfrac 12\right) (pq)^{-2(j+1)} \ell_{j+1}^{-(r_-+(\mu+\rho_-)p+pqb_j+2p+1)} }{ (\mu+\rho_-+qb_j+2)^{2p} (r_-+(\mu+\rho_-)p+pqb_j+2p+1)}  \notag \\ & \qquad \times  (R+t)^{-[r_++(\mu+\rho_+)p+pq a_j]} (t-L_{j+1}T_0)^{r_-+(\mu+\rho_-)p+pqb_j+2p+1} \label{lb F j+1}
\end{align} for $t\in [L_{j+1}T_0,T)$.
Therefore, setting
\begin{align}
B_{j+1} & \doteq \frac{B K^p \left(pq-\tfrac 12\right) (pq)^{-2(j+1)} \ell_{j+1}^{-(r_-+(\mu+\rho_-)p+pqb_j+2p+1)} B_j^{pq} }{ (\mu+\rho_-+qb_j+2)^{2p} (r_-+(\mu+\rho_-)p+pqb_j+2p+1)}  \label{rec rel Bj}, \\
a_{j+1} & \doteq r_++(\mu+\rho_+)p+pq a_j \label{rec rel aj}, \\
b_{j+1} & \doteq r_-+(\mu+\rho_-)p+2p+1+pq b_j \label{rec rel bj},
\end{align} it follows that \eqref{lb F j+1} is just \eqref{lb Fj} for $j+1$.

Let us determine the explicit representations of $a_j$ and $b_j$. Using \eqref{rec rel aj} recursively, we arrive at 
\begin{align*}
a_j & =pq  \, a_{j-1}+ r_++(\mu+\rho_+)p \\ & = (pq)^2  a_{j-2} +(r_++(\mu+\rho_+)p) (1+pq)\\ &= \ldots = (pq)^j a_0+ (r_++(\mu+\rho_+)p) \sum_{k=0}^{j-1}(pq)^k. 
\end{align*} Employing the identity
\begin{align}\label{somma parziale serie geometrica}
\sum_{k=0}^{j-1}(pq)^k = \frac{(pq)^j-1}{pq-1},
\end{align} we conclude that 
\begin{align}
a_j= \left[a_0+\frac{r_++(\mu+\rho_+)p}{pq-1}\right](pq)^j-\frac{r_++(\mu+\rho_+)p}{pq-1} \label{repr aj}.
\end{align} Analogously, from \eqref{rec rel bj} we get
\begin{align}
b_j= \left[b_0+\frac{r_-+(\mu+\rho_-)p+2p+1}{pq-1}\right](pq)^j-\frac{r_-+(\mu+\rho_-)p+2p+1}{pq-1}  \label{repr bj}.
\end{align}
Our next goal is to determine for $j$ large enough a lower bound for $B_j$, whose dependence on $j$ is simpler than the one of $B_j$. We begin by rewriting \eqref{rec rel Bj} as follows:
\begin{align}
B_{j} & = \frac{B K^p \left(pq-\tfrac 12\right) (pq)^{-2j} \ell_{j}^{-b_j} B_{j-1}^{pq} }{ (\mu+\rho_-+qb_{j-1}+2)^{2p} \, b_j}. 
\end{align} By \eqref{rec rel bj}, it follows that 
\begin{align*}
\mu+\rho_-+qb_{j-1} = \frac{1}{p}[b_j-(r_-+2p+1)]\leq \frac{b_j}{p}-2,
\end{align*} which implies
\begin{align*}
(\mu+\rho_-+qb_{j-1}+2)^{2p}\leq p^{-2p}b_j^{2p}.
\end{align*} Consequently, for any $j\in \mathbb{N}$ we showed that
\begin{align}
B_j\geq B K^p \left(pq-\tfrac 12\right) p^{2p} \frac{\ell_j^{-b_j}(pq)^{-2j}}{b_j^{2p+1}}B_{j-1}^{pq}. \label{lb Bj inter}
\end{align} If we denote $E\doteq b_0+\frac{r_-+(\mu+\rho_-)p+2p+1}{pq-1}$, then from \eqref{repr bj} we get
\begin{align}
b_j\leq E (pq)^j \label{ub bj}
\end{align} for any $j\in\mathbb{N}$. Finally, we remark that
\begin{align*}
\lim_{j\to \infty} \ell_{j}^{-b_j}&= \lim_{j\to\infty} \exp(-b_j \ln(\ell _j)) \\
&= \lim_{j\to\infty} \exp\left(\left(-E(pq)^j+(E-b_0)\right) \ln(1+(pq)^{-j})\right) \\
&= \lim_{j\to\infty} \exp\left(-E \,\frac{\ln(1+(pq)^{-j})}{(pq)^{-j}}\right)=\mathrm{e}^{-E}>0.
\end{align*} In particular, there exists $M=M(p,q,\mu,r,\rho,a)>0$ such that
\begin{align}
\ell_{j}^{-b_j}\geq M \label{ub ellj^bj} \quad \mbox{for any} \, j\in\mathbb{N}.
\end{align}
Summarizing, from \eqref{lb Bj inter}, \eqref{ub bj} and \eqref{ub ellj^bj}, we conclude that
\begin{align}
B_j\geq B K^p \left(pq-\tfrac 12\right) p^{2p} M E^{-(2p+1)} (pq)^{-(2p+3)j} B_{j-1}^{pq}= DQ^{-j} B_{j-1}^{pq} \label{lb Bj inter 2}
\end{align} for any $j\in\mathbb{N}$, where $D\doteq B K^p \left(pq-\tfrac 12\right) p^{2p} M E^{-(2p+1)}$ and $Q\doteq (pq)^{2p+3}$.

By applying the logarithmic function to both sides of \eqref{lb Bj inter 2} and using recursively the resulting inequality, we obtain
\begin{align*}
\ln B_j & \geq pq \ln B_{j-1}-j\ln Q+\ln D \notag \\
& \geq (pq)^2 \ln B_{j-2}-[j+(j-1)pq]\ln Q+ (1+pq)\ln D \notag \\
& \geq \ldots  \geq (pq)^j \ln B_{0}-\ln Q\sum_{k=0}^{j-1}(j-k)(pq)^k+\ln D \sum_{k=0}^{j-1}(pq)^k \notag \\
& =\left[\ln A-\frac{pq \ln Q}{(pq-1)^2}+\frac{\ln D}{pq-1}\right] (pq)^j+\frac{\ln Q}{pq-1} \left[\frac{1}{pq-1}+j+1\right]-\frac{\ln D}{pq-1},
\end{align*} where in the last step we used \eqref{somma parziale serie geometrica} and the following identity
\begin{align} \label{identity sum (j-k)(pq)^k}
\sum_{k=0}^{j-1}(j-k)(pq)^k=\frac{1}{pq-1}\left(\frac{(pq)^{j+1}-1}{pq-1}-(j+1)\right).
\end{align}
If we denote by $j_0\in\mathbb{N}$ the smallest nonnegative integer greater than $\frac{\ln D}{\ln Q}-\frac{pq}{pq-1}$ and $H\doteq Q^{-(pq)/(pq-1)^2} D^{1/(pq-1)}$, then for any $j\geq j_0$
\begin{align}
\ln B_j \geq (pq)^j \ln (AH) \label{log ln Bj}.
\end{align}
Our next and final step is to show that for $t$ above a certain $A$-depenent threshold $F(t)$ cannot be finite. To achieve this purpose, we are going to show that the right-hand side of \eqref{lb Fj} blows up as $j\to \infty$. 

Combining \eqref{lb Fj}, \eqref{repr aj}, \eqref{repr bj} and \eqref{log ln Bj}, since $L_j \uparrow L_\infty$ as $j\to \infty$, for $t\in [L_\infty T_0,T)$ and $j\in\mathbb{N}$, $j\geq j_0$ we have
\begin{align}
F(t) & \geq \exp\{(pq)^j \ln(AH)\} (t+R)^{-a_j}(t-L_\infty T_0)^{b_j} \notag\\
    & = \exp\{(pq)^j \ln(AH) -a_j \ln(t+R)+b_j\ln(t-L_\infty T_0)\} \notag \\
    & = \exp\{(pq)^j J(t)\} (t+R)^{\frac{r_++(\mu+\rho_+)p}{pq-1}}(t-L_\infty T_0)^{-\frac{r_-+(\mu+\rho_-)p+2p+1}{pq-1}}, \label{lb F 2nd last}
\end{align} where 
\begin{align*}
J(t) \doteq \ln(AH) -\left(a_0+\frac{r_++(\mu+\rho_+)p}{pq-1}\right)\ln(t+R)+\left(b_0+\frac{r_-+(\mu+\rho_-)p+2p+1}{pq-1}\right)\ln(t-L_\infty T_0).
\end{align*}
If we define $T_1\doteq \max\{R,2 L_\infty T_0\}$, then, for $t\geq T_1$: $\ln(t+R)\leq \ln(2t)$, $\ln(t-L_\infty T_0)\geq \ln\left(\frac{t}{2}\right)$. Hence, for $t\geq T_1$ it results
\begin{align}
J(t) & \geq  \ln(AH) -\left(a_0+\frac{r_++(\mu+\rho_+)p}{pq-1}\right)\ln(2t)+\left(b_0+\frac{r_-+(\mu+\rho_-)p+2p+1}{pq-1}\right)\ln\left(\frac t2\right) \notag\\
&  \geq \ln(AH) +\left(b_0-a_0+\frac{2p+1-\left((r_+-r_-)+(\rho_+-\rho_-)p\right)}{pq-1}\right)\ln t \notag\\
 & \quad -\left(b_0+a_0+\frac{2p+1+(r_++r_-)+(2\mu+\rho_++\rho_-)p}{pq-1}\right) \ln2 \notag \\
& = \ln(AH) + \left(a+\frac{2p+1-(r+\rho p)}{pq-1}\right)\ln t -\left(|a|+\frac{2p+1+|r|+(2\mu+|\rho|) p}{pq-1}\right)\ln 2 \notag \\
& = \ln(A\widetilde{H}) + \left(a+\frac{2p+1-(r+\rho p)}{pq-1}\right)\ln t \doteq \widetilde{J}(t), \label{J>tildeJ}
\end{align} where $\widetilde{H}\doteq 2^{-\big(|a|+\frac{2p+1+|r|+(2\mu+|\rho|) p}{pq-1}\big)}H$.

Let us remark that $\widetilde{J}(t)>0$ if and only if $$t>CA^{-\big(a+\frac{2p+1-(r+\rho p)}{pq-1}\big)^{-1}},$$ where $C\doteq \widetilde{H}^{-\big(a+\frac{2p+1-(r+\rho p)}{pq-1}\big)^{-1}}$.

Since we are working under the assumption $a(pq-1)+2p+1-(\rho p+r)>0$, we can fix $A_0=A_0(a,T_0,R,p,q,r,\rho)>0$ such that for any $A\in(0,A_0]$: $$CA^{-\big(a+\frac{2p+1-(r+\rho p)}{pq-1}\big)^{-1}}\geq \max\{R,2L_\infty T_0\}.$$
Finally, if $A\in(0,A_0]$, $t>CA^{-\big(a+\frac{2p+1-(r+\rho p)}{pq-1}\big)^{-1}}$ and $j\in\mathbb{N}$, $j\geq j_0$, then, $\widetilde{J}(t)>0$ and, due to \eqref{J>tildeJ}, taking the limit as $j\to \infty$ of the right-hand side of \eqref{lb F 2nd last} we conclude that $F(t)$ may not be finite. Therefore, we proved that $F$ blows up in finite time and we derived the upper bound estimat for the lifespan in \eqref{ub est F}. 


\noindent (2) Let us consider the case $\alpha(pq-1)+q+2>rq+\rho$ provided that \eqref{lb G prop} holds. Also in this case \eqref{IF F}-\eqref{IF G} is the iteration frame, and we use it to derive a sequence of lower bound estimates for $G$. The sequence $\{L_j\}_{j\in\mathbb{N}}$ is the same as in (1).

We show that $G$ satisfies the following lower bound estimate for any $j\in\mathbb{N}$:
\begin{align}
G(t)\geq K_j (R+t)^{-\alpha_j} (t-L_jT_0)^{\beta_j} \qquad \mbox{for any} \, t\in[L_jT_0,T), \label{lb Gj}
\end{align} where $\{K_j\}_{j\in\mathbb{N}}$, $\{\alpha_j\}_{j\in\mathbb{N}}$, $\{\beta_j\}_{j\in\mathbb{N}}$ are sequences of nonnegative real numbers that will be determine throughout the iteration argument.

The validity of \eqref{lb Gj} for $j=0$ is guaranteed by \eqref{lb G prop}, setting $K_0\doteq \widetilde{A}$, $a_0\doteq \alpha_-$, $b_0\doteq \alpha_+$. Next we prove the induction step.
Plugging \eqref{lb Gj} into \eqref{IF F}, for $t\in [L_jT_0,T)$ we have
\begin{align*}
F(t) &\geq  B  (R+t)^{-r_+}\, \mathrm{e}^{-t}\int_{L_jT_0}^t\mathrm{e}^{\tau}\int_{L_jT_0}^\tau(s-L_jT_0)^{r_-}(G(s))^p \mathrm{d}s\, \mathrm{d}\tau  \\
&\geq  B K_j^p (R+t)^{-(r_++p\alpha_j)}\, \mathrm{e}^{-t}\int_{L_jT_0}^t\mathrm{e}^{\tau}\int_{L_jT_0}^\tau(s-L_jT_0)^{r_-+p\beta_j} \mathrm{d}s\, \mathrm{d}\tau  \\
&\geq  B K_j^p (r_-+p\beta_j+1)^{-1} (R+t)^{-(r_++p\alpha_j)}\, \mathrm{e}^{-t}\int_{L_jT_0}^t\mathrm{e}^{\tau}(\tau-L_jT_0)^{r_-+p\beta_j+1} \mathrm{d}\tau. 
\end{align*} For $t\in [L_{j+1}T_0,T)$, by shrinking the domain of integration, we may estimate the last $\tau$-integral as follows:
\begin{align*}
& \mathrm{e}^{-t}\int_{L_jT_0}^t\mathrm{e}^{\tau}(\tau-L_jT_0)^{r_-+p\beta_j+1} \mathrm{d}\tau \notag \\
& \quad \geq  \mathrm{e}^{-t}\int_{\tfrac{t}{\ell_{j+1}}}^t\mathrm{e}^{\tau}(\tau-L_jT_0)^{r_-+p\beta_j+1} \mathrm{d}\tau \notag \\
& \quad \geq \ell_{j+1}^{-(r_-+p\beta_j+1)} (t-L_{j+1}T_0)^{r_-+p\beta_j+1} \, \mathrm{e}^{-t}\int_{\tfrac{t}{\ell_{j+1}}}^t\mathrm{e}^{\tau} \mathrm{d}\tau \notag \\
& \quad \geq \ell_{j+1}^{-(r_-+p\beta_j+1)} (t-L_{j+1}T_0)^{r_-+p\beta_j+1} \,  \left( 1-\mathrm{e}^{-\big(1-\frac{1}{\ell_{j+1}}\big)t}\right) \notag \\
& \quad \geq \left(pq-\tfrac{1}{2}\right)\ell_{j+1}^{-(r_-+p\beta_j+1)} (pq)^{-2(j+1)} (t-L_{j+1}T_0)^{r_-+p\beta_j+1},
\end{align*} where we used \eqref{exp factor} in the last inequality.
Summarizing, for $t\in [L_{j+1}T_0,T)$ we proved that
\begin{align}
F(t) &\geq  B \left(pq-\tfrac{1}{2}\right) K_j^p \ell_{j+1}^{-(r_-+p\beta_j+1)} (pq)^{-2(j+1)} (r_-+p\beta_j+1)^{-1} (R+t)^{-(r_++p\alpha_j)}  (t-L_{j+1}T_0)^{r_-+p\beta_j+1}. \label{lb F inter 2}
\end{align}

Using \eqref{lb F inter 2} in \eqref{IF G}, for $t\in[L_{j+1}T_0,T)$ we obtain
\begin{align*}
G(t) & \geq K (R+t)^{-(\mu+\rho_+)} \int_{L_{j+1}T_0}^t\int_{L_{j+1}T_0}^\tau s^{\mu+\rho_-}(F(s))^q \mathrm{d}s\, \mathrm{d}\tau  \\
& \geq K   B^q \left(pq-\tfrac{1}{2}\right)^q K_j^{pq} \ell_{j+1}^{-(r_-+p\beta_j+1)q} (pq)^{-2q(j+1)} (r_-+p\beta_j+1)^{-q}\\ 
& \quad \times (R+t)^{-(\mu+\rho_++(r_++p\alpha_j)q)} \int_{L_{j+1}T_0}^t\int_{L_{j+1}T_0}^\tau (s-L_{j+1}T_0)^{\mu+\rho_-+(r_-+p\beta_j+1)q} \mathrm{d}s\, \mathrm{d}\tau   \\
& \geq \frac{K B^q \left(pq-\tfrac{1}{2}\right)^q K_j^{pq} \ell_{j+1}^{-(r_-+p\beta_j+1)q} (pq)^{-2q(j+1)}}{ (r_-+p\beta_j+1)^{q} (\mu+\rho_-+(r_-+p\beta_j+1)q+2)^2}\\ 
& \quad \times (R+t)^{-(\mu+\rho_++(r_++p\alpha_j)q)} (t-L_{j+1}T_0)^{\mu+\rho_-+(r_-+p\beta_j+1)q+2},
\end{align*} which is exactly \eqref{lb Gj} for $j+1$ provided that 
\begin{align}
K_{j+1} & \doteq \frac{K B^p \left(pq-\tfrac 12\right)^q (pq)^{-2q(j+1)}  \ell_{j+1}^{-(r_-+p\beta_j+1)q} K_j^{pq} }{ (r_-+p\beta_j+1)^{q} (\mu+\rho_-+(r_-+p\beta_j+1)q+2)^2}  \label{rec rel Kj}, \\
\alpha_{j+1} & \doteq \mu+\rho_++r_+q+pq \alpha_j \label{rec rel alphaj}, \\
\beta_{j+1} & \doteq \mu+\rho_-+r_-q+q+2+pq \beta_j \label{rec rel betaj}.
\end{align} 
Employing \eqref{rec rel alphaj} and \eqref{rec rel betaj}, we derive the following representations
\begin{align}
\alpha_j &= \left[\alpha_0+\frac{\mu+\rho_++r_+q}{pq-1}\right](pq)^j-\frac{\mu+\rho_++r_+q}{pq-1}, \label{repr alphaj}\\ 
\beta_j &= \left[\beta_0+\frac{\mu+\rho_-+r_-q+q+2}{pq-1}\right](pq)^j-\frac{\mu+\rho_-+r_-q+q+2}{pq-1}  \label{repr betaj}.
\end{align}
Now we determine a lower bound for $K_j$. By \eqref{rec rel Kj} we get
\begin{align}
K_{j} & = \frac{K B^p \left(pq-\tfrac 12\right)^q }{ (r_-+p\beta_{j-1}+1)^{q} \beta_j^2}  (pq)^{-2qj}  \ell_{j}^{-(r_-+p\beta_{j-1}+1)q} K_{j-1}^{pq} \label{repr Kj}. 
\end{align} From \eqref{rec rel betaj}, we have 
\begin{align*}
r_-+p\beta_{j-1}+1=\frac{1}{q}\left[\beta_j-(\mu+\rho_-+2)\right]\leq \frac{\beta_j}{q}.
\end{align*} Hence, from \eqref{repr Kj} it follows that
\begin{align}
K_j\geq K B^p \left(pq-\tfrac 12\right)^q  q^q (pq)^{-2qj}  \ell_{j}^{-(r_-+p\beta_{j-1}+1)q} \beta_j^{-(q+2)} K_{j-1}^{pq}. \label{lb Kj inter}
\end{align} Let us denote $\widetilde{E}\doteq \beta_0+\frac{\mu+\rho_-+r_-q+q+2}{pq-1}$. By \eqref{repr betaj},
\begin{align}
\beta_j\leq \widetilde{E} (pq)^j \label{ub betaj}
\end{align} for any $j\in\mathbb{N}$. Since
\begin{align*}
\lim_{j\to \infty} \ell_{j}^{-\beta_j}& = \lim_{j\to\infty} \exp\left(\left(-\widetilde{E}(pq)^j+(\widetilde{E}-\beta_0)\right) \ln(1+(pq)^{-j})\right) \\
&= \lim_{j\to\infty} \exp\left(-\widetilde{E} \,\frac{\ln(1+(pq)^{-j})}{(pq)^{-j}}\right)=\mathrm{e}^{-\widetilde{E}}>0,
\end{align*} there exists $\widetilde{M}=\widetilde{M}(p,q,\mu,r,\rho,\alpha)>0$ such that
\begin{align}
\ell_{j}^{-\beta_j}\geq \widetilde{M} \label{ub ellj^betaj} \quad \mbox{for any} \, j\in\mathbb{N}.
\end{align}
Combining \eqref{lb Kj inter}, \eqref{ub betaj} and \eqref{ub ellj^betaj}, we find that
\begin{align}
K_j\geq \widetilde{D} P^{-j}  K_{j-1}^{pq} \label{lb Kj inter 2}
\end{align} for any $j\in\mathbb{N}$, where $\widetilde{D}\doteq K B^p \left(pq-\tfrac 12\right)^q  q^q \widetilde{M} \widetilde{E}^{-(q+2)}$ and $P\doteq (pq)^{3q+2}$.

By \eqref{lb Kj inter 2}, repeating similar computations to those for $\ln B_j$ in the previous case, we obtain
\begin{align*}
\ln K_j & \geq \left[\ln \widetilde{A}-\frac{pq \ln P}{(pq-1)^2}+\frac{\ln \widetilde{D}}{pq-1}\right] (pq)^j+\frac{\ln P}{pq-1} \left[\frac{1}{pq-1}+j+1\right]-\frac{\ln \widetilde{D}}{pq-1}.
\end{align*} Let $j_1\in\mathbb{N}$ be the smallest nonnegative integer greater than $\frac{\ln \widetilde{D}}{\ln P}-\frac{pq}{pq-1}$ and $N\doteq P^{-(pq)/(pq-1)^2} \widetilde{D}^{1/(pq-1)}$. Therefore, for any $j\geq j_1$
\begin{align}
\ln K_j \geq (pq)^j \ln (\widetilde{A}N) \label{log ln Kj}.
\end{align}
Putting together \eqref{lb Gj}, \eqref{repr alphaj}, \eqref{repr betaj} and \eqref{log ln Kj}, for $t\in [L_\infty T_0,T)$ and $j\in\mathbb{N}$, $j\geq j_1$ we have
\begin{align}
G(t) & \geq \exp\{(pq)^j \ln(\widetilde{A}N)\} (t+R)^{-\alpha_j}(t-L_\infty T_0)^{\beta_j} \notag\\
    & = \exp\{(pq)^j I(t)\} (t+R)^{\frac{\mu+\rho_++r_+q}{pq-1}}(t-L_\infty T_0)^{-\frac{\mu+\rho_-+r_-q+q+2}{pq-1}}, \label{lb G 2nd last}
\end{align} where 
\begin{align*}
I(t) \doteq \ln(\widetilde{A}N) -\left(\alpha_0+\frac{\mu+\rho_++r_+q}{pq-1}\right)\ln(t+R)+\left(\beta_0+\frac{\mu+\rho_-+r_-q+q+2}{pq-1}\right)\ln(t-L_\infty T_0).
\end{align*}
For $t\geq T_1= \max\{R,2 L_\infty T_0\}$, it results
\begin{align}
I(t) & \geq \ln(\widetilde{A}N) -\left(\alpha_0+\frac{\mu+\rho_++r_+q}{pq-1}\right)\ln(2t)+\left(\beta_0+\frac{\mu+\rho_-+r_-q+q+2}{pq-1}\right)\ln\left(\frac t2\right) \notag\\
&  \geq \ln(\widetilde{A}N) +\left(\beta_0-\alpha_0+\frac{q+2-\left((\rho_+-\rho_-)+(r_+-r_-)q\right)}{pq-1}\right)\ln t \notag\\
 & \quad - \left(\beta_0+\alpha_0+\frac{q+2+2\mu+(\rho_++\rho_-)+(r_++r_-)q}{pq-1}\right)\ln2 \notag \\
& = \ln(\widetilde{A}N) + \left(\alpha+\frac{q+2-(\rho+r q)}{pq-1}\right)\ln t -\left(|\alpha|+\frac{q+2+2\mu+|\rho|+|r| q}{pq-1}\right)\ln 2 \notag \\
& = \ln(\widetilde{A}\widetilde{N}) + \left(\alpha+\frac{q+2-(\rho+r q)}{pq-1}\right)\ln t \doteq \widetilde{I}(t), \label{I>tildeI}
\end{align} where $\widetilde{N}\doteq 2^{-\big(|\alpha|+\frac{q+2+2\mu+|\rho|+|r| q}{pq-1}\big)}N$.

We point out that $\widetilde{I}(t)>0$ if and only if $$t>\widetilde{C}\widetilde{A}^{-\big(\alpha+\frac{q+2-(\rho+r q)}{pq-1}\big)^{-1}},$$ where $\widetilde{C}\doteq \widetilde{N}^{-\big(\alpha+\frac{q+2-(\rho+r q)}{pq-1}\big)^{-1}}$.

Thanks to the assumption $\alpha(pq-1)+q+2-(rq+\rho)>0$, we can fix $\widetilde{A}_0=\widetilde{A}_0(\alpha,T_0,R,p,q,r,\rho)>0$ such that for any $\widetilde{A}\in(0,\widetilde{A}_0]$: $$\widetilde{C}\widetilde{A}^{-\big(\alpha+\frac{q+2-(\rho+r q)}{pq-1}\big)^{-1}}\geq \max\{R,2L_\infty T_0\}.$$
In the end, taking $\widetilde{A}\in(0,\widetilde{A}_0]$, $t>\widetilde{C}\widetilde{A}^{-\big(\alpha+\frac{q+2-(\rho+r q)}{pq-1}\big)^{-1}}$ and $j\in\mathbb{N}$, $j\geq j_1$, we find that $\widetilde{I}(t)>0$, so letting $j\to \infty$ on the right-hand side of \eqref{lb G 2nd last}, because of \eqref{I>tildeI}, we see that $G(t)$ cannot be finite. Then, we showed that $G$ blows up in finite time and we proved \eqref{ub est G}.


\paragraph*{\textit{Scattering producing case}} As in previous case, we obtain \eqref{IF F} for $F$. By the identity
\begin{align*}
G''(t)+b(t)G'(t) = (m(t))^{-1}\frac{\mathrm{d}}{\mathrm{d}t}\left(m(t)G'(t)\right),
\end{align*} where the multiplier $m$ is defined in \eqref{defm}, the nonegativity of the initial conditions $G(0),G'(0)$ and repeating analogous computations to those in in Section \ref{sce3}, we find
\begin{align}
G(t)&\geq K\int_0^t \int_0^\tau (R+s)^{-\rho}|F(s)|^q \mathrm{d}s\, \mathrm{d}\tau,\label{IF G s.p.}
\end{align} where $K=K(\widetilde{B},b)>0$.
Employing the iteration fram \eqref{IF F}-\eqref{IF G s.p.}, we can repeat the same proof as in the scale-invariant case for $\mu=0$.
\end{proof}

\section{Proof of the main results} \label{Section proof of the thms}

\subsection{Proof of Theorem \ref{theorem1}}\label{proofthem1}

We apply Proposition \ref{Prop comparison argument} when $b(t)=\frac{\mu}{1+t}$ to the pair $(U,V)$ introduced in Subsection \ref{itfr}. Since $U,V$ satisfy the iteration frame \eqref{U1}-\eqref{V1}, by the lower bounds estimates \eqref{1st lb U}, \eqref{2nd lb U}  and \eqref{lb V} we obtain the blow-up range in \eqref{loweralphanmu} and the upper bound estimates for the lifespan in \eqref{ub lifespan thm s.i.}.

\subsection{Proof of Theorem \ref{theorem2}}
In order to prove this theorem, we employ Proposition \ref{Prop comparison argument} when $b\in L^1([0,\infty))$ is a nonnegative function to the pair $(U,V)$ introduced in Subsection \ref{loeboun}. The functionals $U,V$ fulfill the iteration frame \eqref{U1}-\eqref{framev}, therefore, by applying the lower bounds estimates \eqref{lowerU11}, \eqref{lowerV11} and \eqref{2nd lb U}, we get the blow-up range in \eqref{blowup range thm s.p.} and the upper bound estimates for the lifespan in  \eqref{ub lifespan thm s.p.}.

\section{Final remarks and open problem}\label{Section FR&OP}

We have already pointed out that Theorem \ref{theorem1} reduces to the result for the classical Nakao problem obtained in \cite{ChenReissig2021} when $\mu=0$. Similarly, in the scattering producing case $b\in L^1([0,\infty))$, the blow-up range found in Theorem \ref{theorem2} is the same one as for $b=0$.

 On the other hand, for the scale-invariant case, the situation changes when $\mu$ is sufficiently large. From \eqref{loweralphanmu} it follows that for $\mu$ large, the blow-up range does not depend on $\mu$. In this case, we expect that our result in Theorem \ref{theorem1} should not be sharp. Indeed, according to what happens for the semilinear wave equation with scale-invariant damping term, cf. \cite{Ab2015, AbLu2013,Da2021,Wa2014}, we conjecture that there exists $\mu^*(n)>0$ such that for any $\mu\geq \mu^*(n)$ the critical exponent for \eqref{eqs} when $b(t)=\frac{\mu}{1+t}$ is given by the critical curve in \eqref{fpq}.

We stress that, in the present paper, we focused on the blow-up of local in time solutions to \eqref{eqs} when, roughly speaking, the time-dependent coefficient $b$ in the $v$-equation does not modify some asymptotic properties of the solution of the classical Nakao problem \eqref{weak3}. On the other hand, when $b$ is \emph{effective} according to the classification from \cite{Wi07}, the situation is completely different and, for the corresponding single semilinear equation, the critical exponent is no longer related to the Strauss exponent and is instead the Fujita exponent, cf. \cite{DaLuRei2013}. Concerning the case $b$ effective in \eqref{eqs}, some partial results for a weakly coupled system of semilinear damped wave equations with effective and distinct time-dependent coefficients are proved in \cite{Dj2021}.

Furthermore, we point out that, to the best of our knowledge, the global existence of small data solutions to \eqref{weak3} and, more generally, to \eqref{eqs} in the two cases considered in this paper, is an open problem.

Finally, in the forthcoming paper \cite{LiPa2025}, we are going to study the blow-up in finite time for the solutions to weakly coupled systems of Nakao-type with derivative-type nonlinearity, that is, when in \eqref{eqs} the nonlinearity $(|v|^p,|u|^q)$ is replaced by $(|\partial_t v|^p,|\partial_t u|^q)$.

\paragraph*{Acknowledgments}
Y.-Q. Li gratefully acknowledges the valuable opportunity for overseas study provided by the China Scholarship Council (CSC) (Grant No. 202406860050). He also sincerely thanks the Department of Mathematics at the University of Bari, where this paper was completed, for offering a supportive research environment. A. Palmieri is partially supported by the PRIN 2022 project ``Anomalies in partial differential equations and applications'' CUP H53C24000820006. A. Palmieri is member of the \emph{Gruppo Nazionale per L'Analisi Matematica, la Probabilit\`{a} e le loro Applicazioni} (GNAMPA) of the \emph{Instituto Nazionale di Alta Matematica} (INdAM).

\end{document}